\documentclass[a4paper]{article}
\usepackage{amsthm}
\usepackage{natbib}
\usepackage{hyperref}
\usepackage{enumitem}
\usepackage[small]{titlesec}
\usepackage{graphicx}

\def \A{{\rm A}}
\def \B{{\rm B}}
\def \O{{\rm O}}
\def \F{{\rm F}}
\def \G{{\rm G}}
\def \Q{{\rm Q}}
\def \MU{M}
\def \NU{N}

\def\R{{{\rm I}\!{\rm R}}}
\def\lint{\mathop{\scriptstyle{\rfloor}}}
    
\defcitealias{gauss}{book}

\defcitealias{lagrange3}{first}

\title{\large{\bf Lambert's Theorem: Geometry or Dynamics?}}

\author{\normalsize{\bf Alain Albouy}\\
\\
\normalsize{IMCCE, UMR 8028,}\\
\normalsize{77, avenue Denfert-Rochereau}\\
\normalsize{75014 Paris, France}\\
\normalsize{Alain.Albouy@obspm.fr}\\
\\
\\
\normalsize{Dedicated to Christian Marchal for his 80$^{\rm th}$ birthday}}

\date{}

\begin{document}

\pagestyle{empty}

\setcounter{page}{0} 

\centerline{\bf Lambert's Theorem: Geometry or Dynamics?}\bigskip

\centerline{Alain Albouy}
\centerline{IMCCE, Observatoire de Paris, UMR 8028, CNRS}
\centerline{77, avenue Denfert-Rochereau, 75014 Paris, France}
\centerline{Alain.Albouy@obspm.fr}\bigskip\bigskip

\bigskip

\noindent This is a post-peer-review, pre-copyedit version of an article published in {\it Celestial Mechanics and Dynamical Astronomy}, as part of the topical collection {\it 50 years of Celestial Mechanics and Dynamical Astronomy}. The final authenticated version is available online at \href{http://dx.doi.org/10.1007/s10569-019-9916-2}{http://dx.doi.org/10.1007/s10569-019-9916-2}. This article should be cited as:

\bigskip
\noindent A.\ Albouy, Lambert's Theorem: Geometry or Dynamics?, {\it Celestial Mechanics and Dynamical Astronomy}, 131 (2019), 40

\tableofcontents

\maketitle

\section*{Abstract}
Lambert's theorem (1761) on the elapsed time along a Keplerian arc drew  the attention of several prestigious mathematicians. In particular, they tried to give simple
and transparent proofs of it (see our timeline Sect.\thinspace\ref{sect9}). We give two new proofs. The first one (Sect.\thinspace\ref{sect4}) goes along the lines of Hamilton's variational proof in his famous paper of 1834, but we shorten his computation in such a way that the hypothesis is now used without redundancy. The second (Sect.\thinspace\ref{sect6}) is among the few which are close to Lambert's geometrical proof. It starts with the new remark that two Keplerian arcs related by the hypothesis of Lambert's theorem correspond to each other through an affine map. We also show (Sect.\thinspace\ref{sect7}) that despite the singularities due to the occurrence of collisions, the classes of arcs related by  Lambert's theorem all have the same topology. We give (Sect.\thinspace\ref{sect8})  some simple related results about conic sections and affine transformations.

\bigskip

\section{Preliminaries}
\label{sect1}
Lambert's theorem is about {\it arcs} of {\it Keplerian orbits}, which we will also call {\it Keplerian arcs}. The statement of this old theorem is more uniform if we understand that Keplerian orbits are {\it extended beyond collisions}. This classical extension is not smooth and is often explained by introducing a regularization. Here we describe it directly in the simplest possible way. We take this opportunity to introduce some basic definitions. The reader can skip this section for a quicker access to the statement of the theorem.

\medskip\noindent{\bf Definition 1.} Call $\O$ the origin of the Euclidean vector space $\R^d$. An {\it extended solution} of Newton's differential system
\begin{equation}\label{n1}
{d^2 q\over dt^2} =-{q\over r^3},\quad \hbox{where } r=\|q\|,
\end{equation}
is a continuous path $\R\to \R^d$, $t\mapsto q$, such that $q=\O$ on a discrete subset of $\R$, and which is an analytic solution of (\ref{n1}) outside of this subset. If $q=\O$ at a time $t_\O$, the extension is characterized as follows: the position $q$ remains on the same ray, and the energy $H$ takes the same value for all $t$ such that $q\neq \O$.
\medskip

{\bf Terminology and notation.} The {\it velocity vector}, first derivative of the {\it position} $q$ with respect to the {\it time} $t$, is denoted by $\dot q$ or $dq/dt$ or $v$. The  {\it acceleration vector} is $\ddot q=d^2q/dt^2$.  The {\it energy} 
\begin{equation}\label{n2}
H={1\over 2}\|\dot q\|^2-{1\over r},
\end{equation}
 is constant along the solutions of (\ref{n1}). A {\it ray} is a closed half-line extending from the origin $\O$. An extended solution which remains on a ray is called a {\it rectilinear solution}.

What happens in a rectilinear solution, for the isolated values of $t$ when $q=\O$, is called a {\it collision} of $q$ with $\O$. Then $q$ is  ``bouncing'' off of $\O$. Explicitly, for such a time $t_\O$, $q=\O$ and the velocity $\dot q$ is infinite. For $t<t_\O$  and $t_\O-t$ sufficiently small, we are on a solution of (\ref{n1}) where the vector $\dot q$ points toward $\O$. It tends to infinity when $t\to t_\O$. For $t_\O<t$ and $t-t_\O$ sufficiently small, we are on a solution of (\ref{n1}) where $\dot q$ points in the opposite direction. It tends to infinity when $t\to t_\O$. 

If, at a given time, the  position vector $q$ and the velocity $\dot q$ span a two-dimensional space, then the solutions of (\ref{n1}) are defined for all time. Such maximal solutions are by definition the extended solutions. As is well known, these are planar solutions which are, if $H<0$, the elliptic solutions described by Kepler, if $H=0$, parabolic solutions, and if $H>0$, hyperbolic solutions.

 \vspace{0.5cm}
\centerline{\includegraphics[width=60mm]{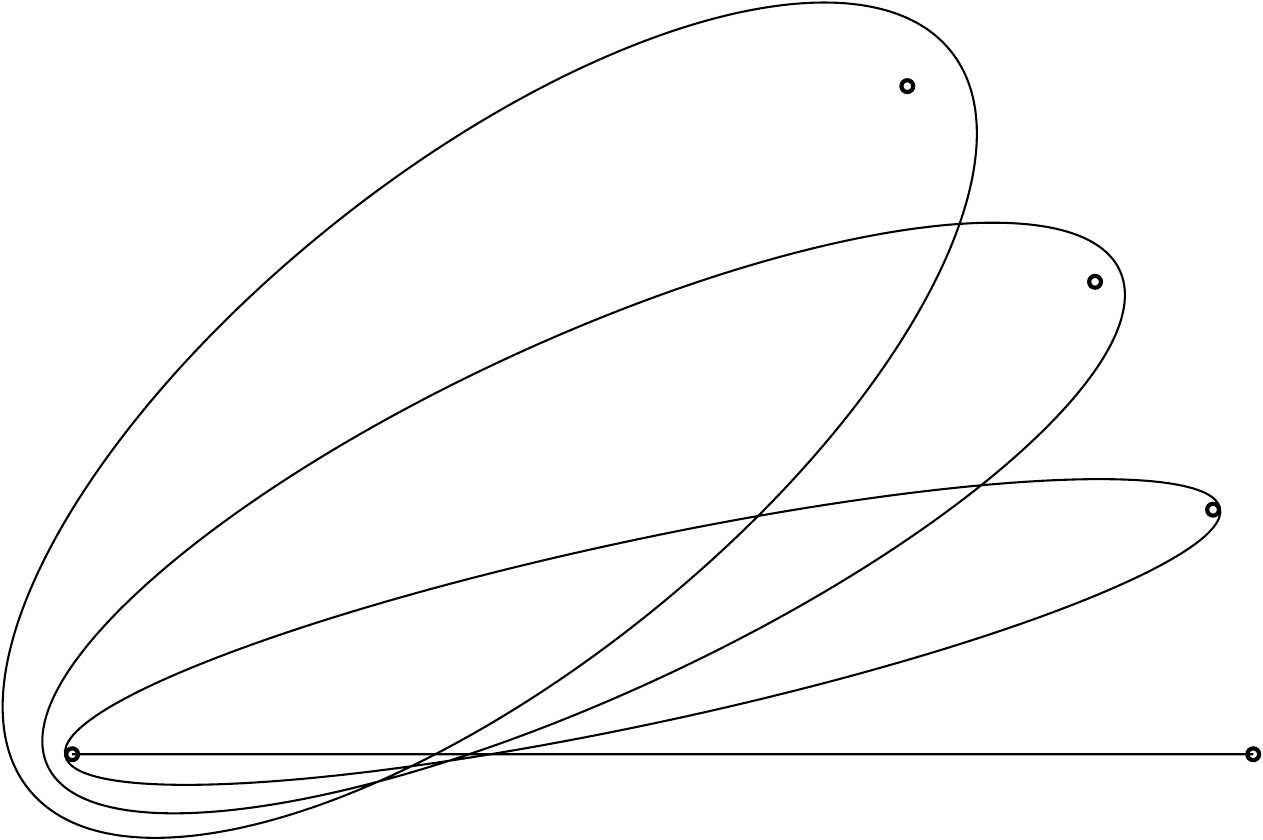}}
\centerline{\sl Fig.\thinspace1. Four Keplerian ellipses with same energy.}
\vspace{0.5cm}

We may justify the extension beyond collision as follows. Consider a one-parameter family of usual solutions  which tends to a rectilinear solution in any way, for example as in Fig.\thinspace1. In the limit the behavior is as we just described.

Figure 2 displays the function $t\mapsto r$ for a bounded solution. What is drawn is the path $u\mapsto (t,r)=(u-\sin u,1-\cos u)$. We recognize a cycloid. The variable $u$ is the {\it eccentric anomaly}, and $t=u-\sin u$ is the {\it Kepler equation} in the rectilinear case. Note that $H=-1/2$ along this solution.

System (\ref{n1}) is autonomous: if $t\mapsto q(t)$ is a solution, then $t\mapsto q(t+\tau)$, where $\tau\in\R$, is also a solution. All the solutions obtained from each other by such a time shift form a class that we call an {\it orbit} or a {\it trajectory}.

\medskip\noindent{\bf Definition 2.} A {\it Keplerian orbit} around $\O$ is a class formed by an extended solution of (\ref{n1}) and all the solutions obtained from it by a time shift. An {\it arc} of Keplerian orbit around $\O$ is a class formed by an extended solution {\it restricted to an interval} $[t_\A,t_\B]$, and all the restricted solutions obtained from it by a time shift. Here $t_\A\in\R$ is called the initial time and $t_\B>t_\A$ the final time. The {\it elapsed time} $\Delta t=t_\B-t_\A$ is invariant by a time shift. The {\it ends} of an arc are the initial position $\A\in\R^d$ and the final position $\B\in\R^d$.  We do not assume that $\A$, $\B$ and $\O$ are distinct.
\medskip

\vspace{0.5cm}
\centerline{\includegraphics[width=30mm]{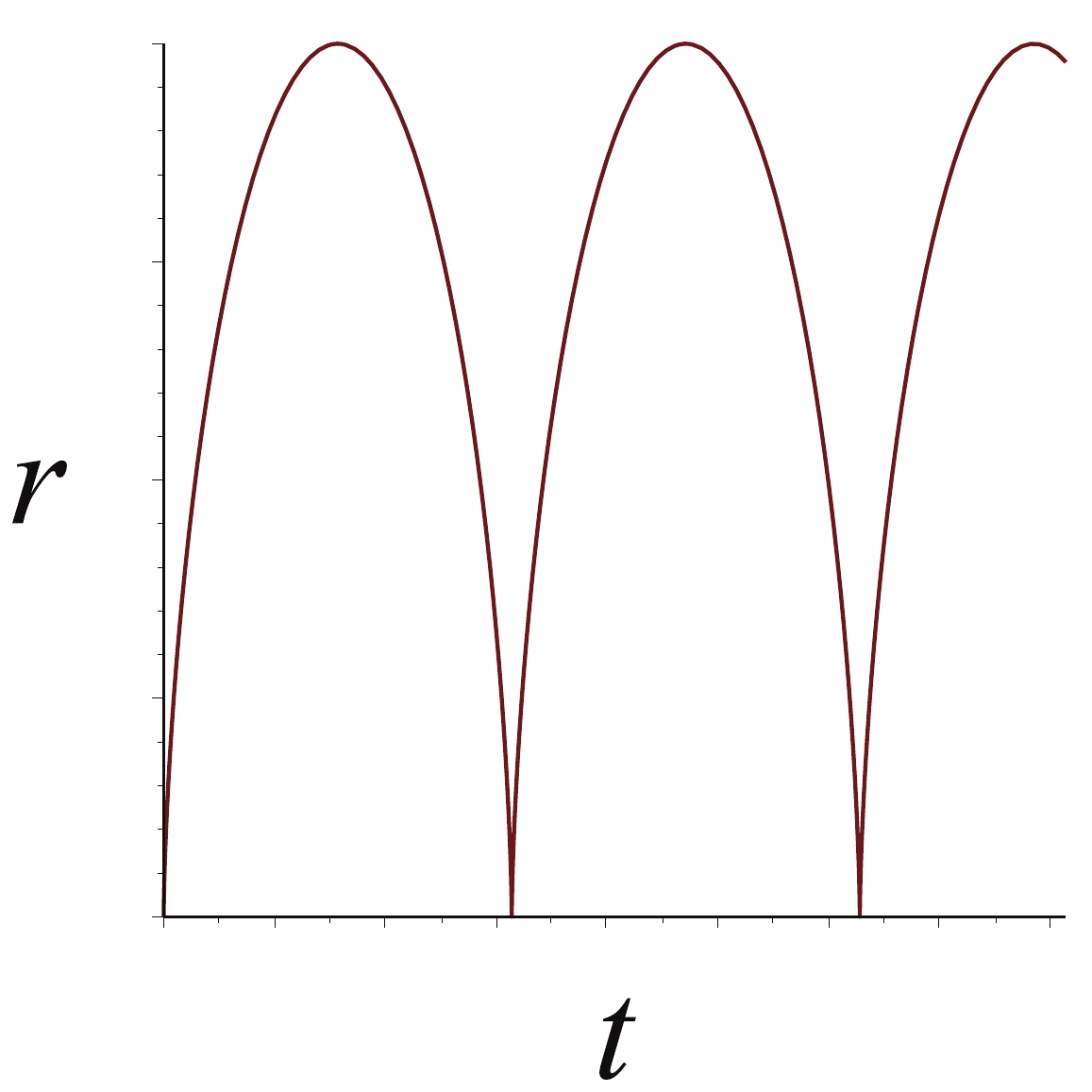}}
\centerline{\sl Fig.\thinspace2. A rectilinear solution with negative energy.}
\vspace{0.5cm}

A Keplerian orbit is a dynamical object: a point $q$ moves according to law (\ref{n1}). We are sometimes only interested in the planar curve described by the body. 

\medskip\noindent{\bf Definition 3.} We call a {\it Keplerian branch} around $\O$ the image of the map $\R\to\R^d, t\mapsto q$ in an extended solution of (\ref{n1}).
\medskip

A Keplerian branch is an ellipse, a parabola, a branch of hyperbola, a compact interval or a ray. The origin $\O$ is always in the convex hull of the branch. When the orbit is rectilinear, the conic section is a double line, and $\O$ is an end of the interval. A conic section in the plane which is not a pair of lines is called {\it irreducible}. Here is a well-known statement.

 \vspace{0.5cm}
\centerline{\includegraphics[width=120mm]{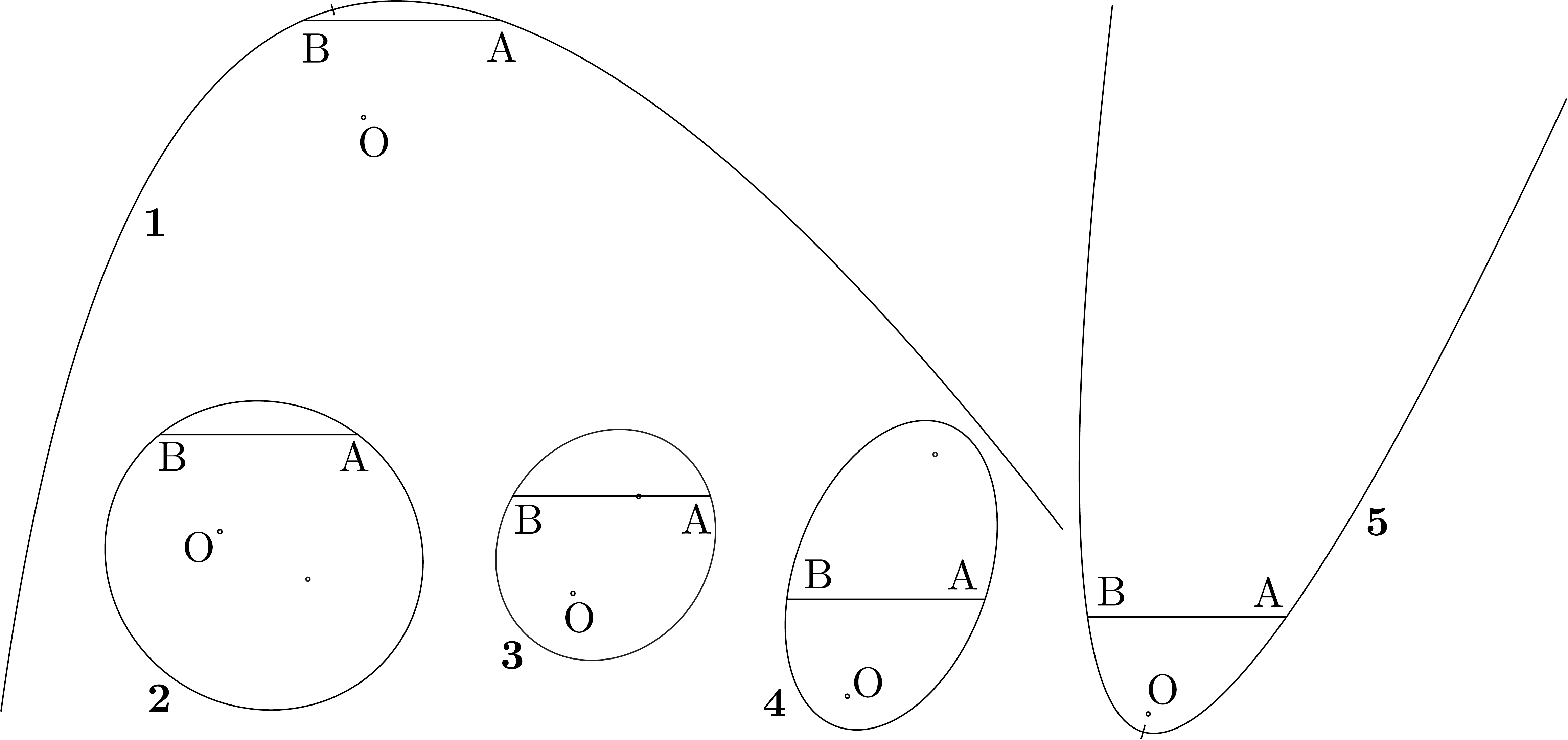}}
\centerline{\sl Fig.\thinspace3. A sequence of Keplerian branches passing through fixed $\A$ and $\B$.}
\vspace{0.5cm}

\medskip\noindent{\bf Lemma 1.} Any irreducible conic section in a Euclidean plane, with a focus at $\O$, includes a unique Keplerian branch around $\O$. This branch is described by two and only two Keplerian orbits around $\O$, which are nonrectilinear and differ only in their orientation. 
\medskip

A natural projection of the {\it space of Keplerian arcs} ${\cal A}$ onto the space of pairs $(\A, \B)\in \R^d\times \R^d$ associates  to any arc its ends.
Figure 3 displays five elements of the family of Keplerian branches passing though two given points $\A$ and $\B$, 1 and 5 being parabolas, 2, 3, 4 being ellipses. We did not represent the branches of hyperbolas which continue the family before the parabola 1 and after the parabola 5. We can follow the one-parameter family of upper arcs, going from $\A$ to $\B$ counterclockwise, from 1 to 4. In the parabola 5, this arc just disappeared at infinity. We can follow the one-parameter family of lower arcs, going from $\A$ to $\B$ clockwise, from 2 to 5.

 \medskip\noindent{\bf Lemma 2.} If the ends $\A$ and $\B$ of an arc are distinct and on the same ray, then the arc is rectilinear.
 
 \medskip\noindent{\bf Proof.} If a nonrectilinear Keplerian orbit crosses a ray again, it crosses it at the same point (with the same velocity): it is an elliptic orbit.\qed

\section{The statement}
\label{sect2}
Lambert's theorem may be stated as a property of some families of Keplerian arcs (see Fig.\thinspace4). We may think of these families as paths in the space ${\cal A}$ of Keplerian arcs. To such a path ${\cal I}\to {\cal A}$, $s\mapsto {\Gamma\!}_s$, where ${\cal I}$ is an open interval, is associated a path ${\cal I}\to \R^d\times \R^d$, $s\mapsto (\A_s,\B_s)$, which describes the displacement of the ends of the arc ${\Gamma\!}_s$.
We may also avoid such technical words by naming a path $s\mapsto {\Gamma\!}_s$ a {\it continuous change} of arc.

\vspace{0.5cm}
\centerline{\includegraphics[width=65mm]{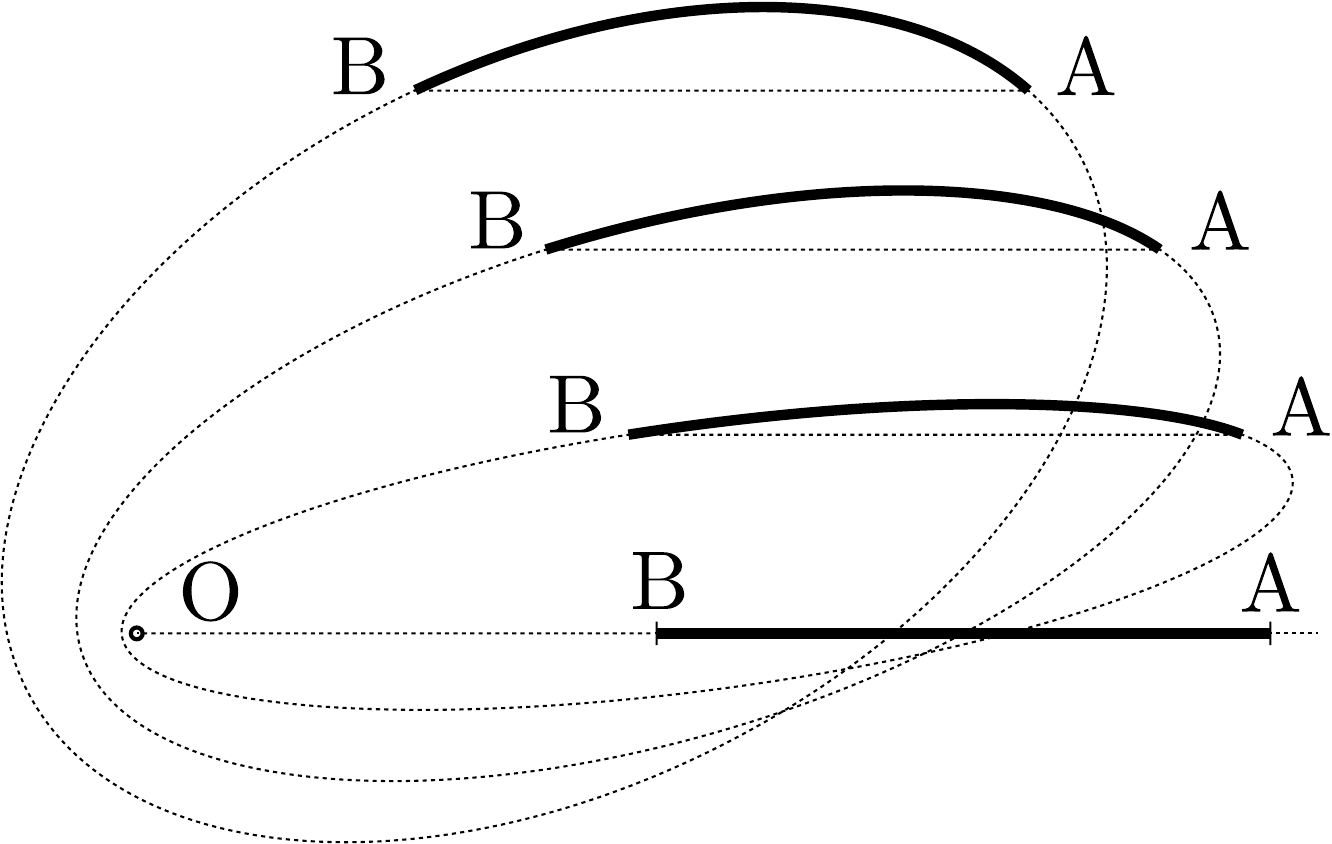}}
\centerline{\sl Fig.\thinspace4. Four Keplerian arcs with same $\|\A\B\|$, $\|\O\A\|+\|\O\B\|$ and $H$.}
\vspace{0.5cm}

\medskip\noindent{\bf Theorem 1 (Lambert).} Consider the Keplerian arcs around the origin $\O$ of $\R^d$. If we change continuously such an arc while keeping constant the distance $\|\A\B\|$ between both ends, the sum of the radii $\|\O\A\|+\|\O\B\|$ and the energy $H$, then the elapsed time $\Delta t=t_\B-t_\A$ is also constant.

\medskip\noindent{\bf Theorem 2 (Lambert).} Starting from any given Keplerian arc, we can arrive at some rectilinear arc by a continuous change which keeps constant the three quantities $\|\A\B\|$, $\|\O\A\|+\|\O\B\|$ and $H$.

\medskip\noindent{\bf Remark 1.} Theorem 2 is usually absent from the statements classically called Lambert's theorem. But the classical authors explain, after crediting to Lambert a statement similar to Theorem 1, how to deduce elegant and useful formulas for $\Delta t$, by reducing the general case to the rectilinear case. They use Theorem~2 but often neglect its proof.

\medskip\noindent{\bf Remark 2.} The classical statement of Theorem 1 is: {\it $\Delta t$ is a function of $H$, $\|\A\B\|$ and $\|\O\A\|+\|\O\B\|$}. However, the ``function'' is ramified and multivalued for two reasons.
Firstly, in the case of an ellipse, $q$ can go from $\A$ to  $\B$ clockwise or counterclockwise, and can make several turns, which gives various arcs with the same $H$ and different $\Delta t$.
Secondly, if $\A$ and $\B$ are not on the same ray and if $H>H_{\min}$, with 
\begin{equation}\label{Hmin}
H_{\min}=-{2\over \|\A\B\|+\|\O\A\|+\|\O\B\|},
\end{equation}
then there are exactly two distinct Keplerian branches with energy $H$ passing through $\A$ and $\B$. 

\medskip\noindent{\bf Remark 3.} A construction indicated as a footnote in \citet{gauss}, \S 106, explains the latter bivaluation in the case of elliptic orbits. Let us recall the classical formula for the energy: 
\begin{equation}\label{Ha}
H=-{1\over 2a},
\end{equation} 
where $a$ is the semimajor axis of the ellipse. If we know a focus and the semimajor axis, the ellipse is characterized by the second focus $\F$, and drawn with a pencil and a string of length $2a$ attached to both foci. Let us draw a circle of center $\A$ and of radius $2a-\|\O\A\|$, and a  circle of center $\B$ and of radius $2a-\|\O\B\|$ (see Fig.\thinspace5). Each intersection is a possible $\F$. The corresponding ellipses have the same semimajor axis $a$ and the same energy $H=-(2a)^{-1}$. 

\vspace{0.5cm}
\centerline{\includegraphics[width=30mm]{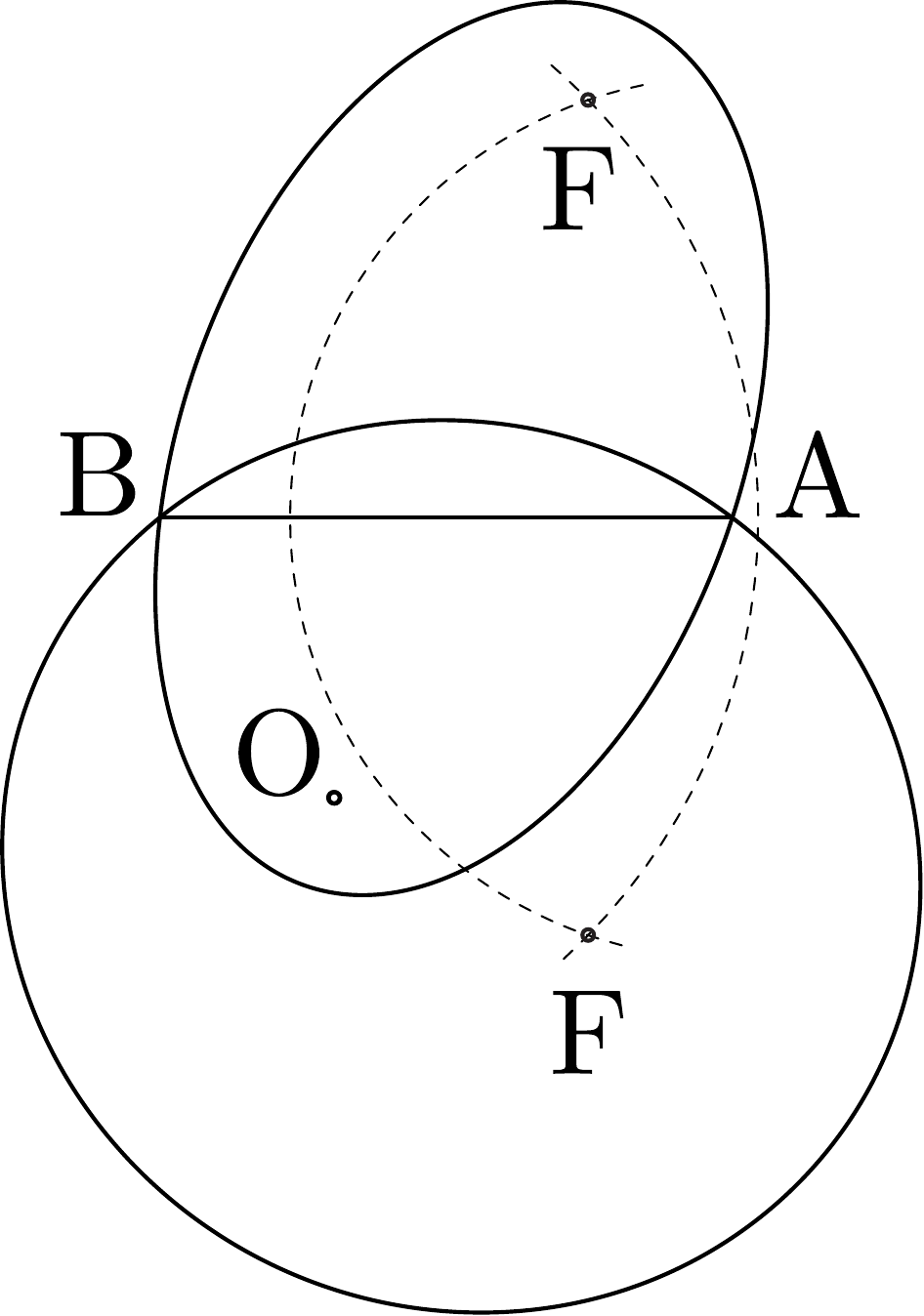}}
\centerline{\sl Fig.\thinspace5. Two Keplerian branches with same energy.}
\vspace{0.5cm}

The circles are tangent when $\pm(2a-\|\O\A\|)=\|\A\B\|\pm(2a-\|\O\B\|)$. The case $(-,+)$ is forbidden since the quantities inside the parentheses are positive. There remains  $\|\O\A\|+\|\A\B\|=\|\O\B\|$, $\|\O\B\|+\|\A\B\|=\|\O\A\|$ or
$4a= \|\O\A\|+\|\O\B\|+\|\A\B\|$. The first cases give rectilinear orbits, which we excluded. The last case corresponds to an $\F$ on the chord $\A\B$ and to the minimal energy $H_{\min}$ in  (\ref{Hmin}). The circles have two distinct intersections if and only if $4a> \|\O\A\|+\|\O\B\|+\|\A\B\|$. In Fig.\thinspace3, ellipse 3 has minimal energy, and ellipses 2 and 4 have the same energy.

\medskip\noindent{\bf Remark 4.} The geometry of the minimal energy described in the previous remark is well known in space dynamics (see \citealt{battin}, figure 3.1). But it may also be checked by throwing metal balls in a room, or in any situation where the nongravitational forces may be neglected. One may consider the trajectories as parabolic, or as elliptic with a focus at the center of the Earth. {\it If we throw such a ball from a point $\A$, and manage to reach a point $\B$  by starting with the least possible velocity, then the focus of the parabolic trajectory is on the chord $\A\B$, whatever the respective altitudes of $\A$ and $\B$.}
This property of the focus of the parabola could be taught in elementary classes. But it appears to be generally ignored.

\medskip\noindent{\bf Remark 5.} We can restate Theorem 1 as follows: {\it the four functions $\|\A\B\|$, $\|\O\A\|+\|\O\B\|$, $H$ and $\Delta t$ are functionally dependent on the space of Keplerian arcs}.  This suggests in turn another statement: {\it $H$ is a function of $\|\A\B\|$, $\|\O\A\|+\|\O\B\|$ and $\Delta t$}. The  second multivaluation mentioned in remark 2 disappears, the elapsed time $\Delta t$ being a uniform parameter for the two families of arcs with given ends described at the end of Sect.\thinspace\ref{sect1}.
We will recall in Sect.\thinspace\ref{sect4} that another parameter can replace $\Delta t$ or $H$ in the statements: the Maupertuis action $w$ of the arc.

\medskip\noindent{\bf Remark 6.} That $\Delta t$ is a uniform parameter for the arcs with given ends is a result of \citet{simo}: {\it In a plane with origin $\O$, for any $\A$ and $\B$ not on the same ray from $\O$, and any $\Delta t>0$, there are exactly two Keplerian arcs around $\O$ going from $\A$ to $\B$ in a time $\Delta t$ and in less than one turn, one clockwise and the other counterclockwise.}

\medskip\noindent{\bf Remark 7.} The above problem of finding a Keplerian arc with given ends and given $\Delta t$ began to be called {\it the Lambert problem} in the 1960's. It was considered earlier as a step in some methods of orbit determination from three or more observations. It found new applications during the conquest of space, and received a considerable attention. What are the relations between the Lambert problem and Lambert's theorem? Lambert's theorem allows to reduce the Lambert problem for a general triangle $\O\A\B$ to the case where $\O\A\B$ are collinear in this order. This reduction is of course optional, and direct methods may be preferred. \citet{gauss} proposes a method of orbit determination where the Lambert problem is a step. He recommends a method to solve the Lambert problem which he begins to explain in his \S 88. He is not very explicit about the relation of his method with Lambert's theorem. \citet{battin} has an influential chapter (which we cited in remark 4) where Lambert's theorem is explained and related to the Lambert problem. Methods are discussed, including Gauss's method, which is related to Lagrange-Gauss-Adams's proof of Lambert's theorem. The above result by Sim\'o is anticipated in a single sentence: ``However, when the time of flight is also given, then, in general, the orbit will be unique.''

\medskip\noindent{\bf Remark 8.} The Lambert problem was indeed posed by Lambert. We can read in \citet{bopp}, p.\ 24, statements that Lambert sent to Euler about the possible sets of data that can be used to determine uniquely a Keplerian orbit:

\noindent``J'ai oubli{\'e} de tourner le probleme \S 210 c'est que l'orbite se trouvera

\noindent 1$^o$ par les 3 cot{\'e}s $FN$, $FM$, $NM$ et le tems $T$ emploi{\'e} {\`a} parcourir l'arc $NM$.

\noindent 2$^o$ par le rapport $(FM:FN)$, l'angle $NFM$, le tems $T$, et le tems periodique.

\noindent Si le diametre du Soleil peut {\^e}tre mesur{\'e} assez exactement, deux observations suffisent pour d{\'e}terminer l'orbite de la Terre par ce dernier th{\'e}oreme\footnote{I forgot to turn problem \S 210 around, that the orbit may be found

\noindent 1$^o$ by the 3 sides $FN$, $FM$, $NM$ and the time $T$ required to traverse the arc $NM$.

\noindent 2$^o$ by the ratio $(FM:FN)$, the angle $NFM$, the time $T$, and the periodic time.

\noindent If the diameter of the Sun can be measured precisely enough, two observations suffice to determine the Earth's orbit using the latter theorem.}.''

\section{The eccentricity vector and the unifocal equation} 
\label{sect3}
There are many ways to prove that the Keplerian orbits are solutions of Newton's differential system (\ref{n1}).  The following known method is extremely brief. We consider the orbits in a plane $\O xy$, write the vector $q=(x,y)$, and set
\begin{equation}\label{o8}
C=x\dot y-y\dot x,\quad \alpha={x\over r}-\dot y C,\quad \beta={y\over r}+\dot x C,\quad E=(\alpha,\beta).
\end{equation}
We check that  $\dot C=0$ and $\dot E=0$. So, the {\it angular momentum} or {\it areal constant} $C$ is constant along the solutions of (\ref{n1}), and the same is true of the {\it eccentricity vector} $E$.  The norm of $E$ is the eccentricity. Its direction is always opposite to the pericenter. To see this, we deduce from (\ref{o8}) that
$$\alpha x+\beta y=r-C^2.$$
Setting $\gamma=C^2$, this is
\begin{equation}\label{unifoc}
r=\alpha x +\beta y+\gamma.
\end{equation}
According to the famous focus-directrix description of a conic section, when the semiparameter $\gamma>0$, this is a branch of conic section with focus at the origin $\O$. The directrix is the line $0 = \alpha x + \beta y + \gamma$. The right-hand side is the distance to the directrix multiplied by the eccentricity $\sqrt{\alpha^2+\beta^2}$. The left-hand side is the distance to $\O$. Pappus proved that a curve described in this way is a conic section in his report about the {\it Surface-loci}, a lost book by Euclid (see \citealt{thomas}, p.\ 493, \citealt{heath}, p.\ 243, \citealt{chasles}, p.\ 44). 
The branch is the whole conic section, except in the case of the hyperbola, where it is the branch whose convex hull contains $\O$.

So, any curve drawn by a nonrectilinear solution of (\ref{n1}) satisfies Eq.\thinspace(\ref{unifoc}). The energy (\ref{n2}) is related to the eccentricity and the angular momentum by
\begin{equation}\label{H}
\alpha^2+\beta^2-1=-{2 C^2\over r}+(\dot x^2+\dot y^2)C^2=2H\gamma.
\end{equation}
Gauss's opinion about Eq.\thinspace(\ref{unifoc}), which we call the {\it unifocal equation} of conic sections, appears in his \citetalias{gauss}, \S 3, p.\ 3:

``Inquiries into the motions of the heavenly bodies, so far as they take place in conic sections, by no means demand a complete theory of this class of curves; but a single general equation rather, on which all others can be based, will answer our purpose. And it appears to be particularly advantageous to select that one to which, while investigating the curve described according to the law of attraction, we are conducted as a characteristic equation. [...] if we denote the distance of the body from the sun by $r$ (always positive), we shall have between $r$, $x$, $y$, the linear equation $r+\alpha x+\beta y=\gamma$, in which $\alpha$, $\beta$, $\gamma$ represent constant quantities, $\gamma$ being from the nature of the case always positive.''

\medskip\noindent{\bf Remark 9.} The signs in our Eq.\thinspace(\ref{unifoc}) differ from Gauss's. We take them from \citet{lagrange5}, \S 7, who presented the above integration of the Kepler problem, which consists in deducing (\ref{unifoc}) from expression (\ref{o8}) of the eccentricity vector. Lagrange is not exactly the first to present this deduction. Jacob Herman published a famous note in 1710, where, however, he did not write $(\alpha,\beta)$ but only the equation $\beta=0$ after a choice of axis (see \citealt{albouy2} for an explanation).

\medskip\noindent{\bf Remark 10.} Several vectorial systems of notation were used to write $(\ref{o8})$. The three-dimensional case was presented with cross products in \citet{gibbs}, \S 61:
$$\quad C=q\times v, \quad E={q\over r}-v\times C.$$
\citet{CushmanDuistermaat} used exterior algebra with vectors and covectors identified through the Euclidean form. This reads
\begin{equation}\label{old5}
C=q\wedge v, \quad E={q\over r}+v\lint C,
\end{equation}
where $C$ is a bivector and $\lint$ is the contracted product. 
We have $v\lint (q\wedge v)=\langle v,q\rangle v-\langle v,v\rangle q$.
This notation is valid in any dimension.

\medskip\noindent{\bf Remark 11.} If we divide (\ref{unifoc}) by $r$ and set $x=r\cos\theta$, $y=r\sin\theta$, $\alpha=-e\cos\theta_0$, $\beta=-e\sin\theta_0$, we find the famous polar equation of Keplerian branches
$${\gamma\over r}=1+e\cos (\theta-\theta_0).$$ This equation in itself is as good as (\ref{unifoc}) for the description of a Keplerian motion, but its form suggests  parameters other than $(\alpha,\beta,\gamma)$. This triple of ``affine parameters'' is a key to all the simplifications reported in the present work. The simplicity of the following Lemma is an illustration.

\medskip\noindent{\bf Lemma 3.} There is a one-to-one map from the space of irreducible planar Keplerian branches (Definition 3) around $\O$ onto the open half-space ${\cal H}=\R\times \R\times ]0,+\infty[$, which consists in associating to a branch the parameters $(\alpha,\beta,\gamma)$ of its unifocal Eq.\thinspace(\ref{unifoc}).  The space of irreducible Keplerian branches which pass through a point $\A\neq \O$ is the intersection with ${\cal H}$ of a plane which cuts the boundary of ${\cal H}$. The space of irreducible Keplerian branches which pass through two points $\A$ and $\B$, not located on the same ray, is the intersection with ${\cal H}$ of a straight line. This line cuts the boundary of ${\cal H}$, except if $\O$ is on the interval $]\A,\B[$, in which case the line is in ${\cal H}$ and $\gamma=C^2$ is constant on the line.

\medskip\noindent{\bf Proof.} If  $(\alpha,\beta,\gamma)$ is given, then we choose $(x,y)$ on the curve and compute $(\dot x,\dot y)$ by making $C=\sqrt \gamma$ in (\ref{o8}). The three conditions are consistent. We get an initial condition for  (\ref{n1}) and a Keplerian branch. The remaining statements are deduced from the equations $r_\A=\alpha x_\A+\beta y_\A+\gamma$ and $r_\B=\alpha x_\B+\beta y_\B+\gamma$.\qed

\section{A minimal proof of Theorem 1}
\label{sect4}
In a famous work inspired by optics, \citet{hamilton1} presented a new approach to dynamics, and wished to demonstrate its effectiveness by a new presentation of Lambert's theorem. He actually deduces new properties and relations which are probably the most original results on the subject after Lambert's book. But then he proposes a new proof which is somewhat disappointing: being similar to a well-known argument introduced by Lagrange, Hamilton's argument does not show any significant simplification. We will present a proof based on the same general variational formula (\ref{3n}), but treating the Kepler problem in a much shorter way.

{\bf Hamilton's variational idea for proving Lambert's theorem.} Denote by $v=\dot q$ the velocity and consider 
\begin{equation}\label{w}
w=\int_{t_\A}^{t_\B} \|v\|^2dt,
\end{equation}
which is stationary on any solution of (\ref{n1}), considered among the paths $[t_\A,t_\B]\to \R^d$, $t\mapsto q$ with arbitrary values  $t_\A$ and $t_\B$, but with same ends $\A$ and $\B$ and same energy $H$.

What we just stated is a variational principle for a natural system, called the Maupertuis principle. The {\it action integral} (\ref{w}) is sometimes called the Maupertuis action. Hamilton called  it the {\it characteristic function} and proved the formula:
\begin{equation}\label{3n}
\delta w=\langle \delta \B,v_\B\rangle-\langle \delta \A,v_\A\rangle+(t_\B-t_\A)\delta H.
\end{equation}
Formula (\ref{3n}) is true for a {\it variation} among the {\it solutions} of the equation of motion (\ref{n1}). A Keplerian arc whose ends are called $\A$ and $\B$ is embedded into a one-parameter family of Keplerian arcs, with varying ends and energy.  In these traditional variation formulas, the parameter of the family is called a variation parameter, and  $\delta f$, the ``infinitesimal variation'' of a function $f$, is the first derivative  of  $f$ with respect to this parameter. Thus, $\delta w$ and $\delta H$ are two numbers, and $\delta \A$ and $\delta \B$ are two vectors tangent to the configuration space, respectively at $\A$ and at $\B$. The velocities $v_\A$ and $v_\B$ of the body at $\A$ and $\B$ are also tangent vectors. The energy $H$ is constant along the solutions, but may vary during the variation.

Formula (\ref{3n})  presents good features for a proof of Lambert's theorem. The variations allowed are exactly what we called the continuous changes. Moreover, the formula is  extremely simple when the energy is constant during the variation. We will use (\ref{3n}) as Hamilton did, but our proof of Proposition 1 is much shorter and treats in a single case the three kinds of conic sections.

\medskip\noindent{\bf Proposition 1 (Hamilton).} If we change continuously a Keplerian arc while keeping constant  $\|\A\B\|$, $\|\O\A\|+\|\O\B\|$ and $H$, then the Maupertuis action $w$ is also constant.

\medskip\noindent{\bf Proof.}  We will prove that the infinitesimal variation $\delta w$ vanishes. As $H$ is fixed, $(\ref{3n})$ becomes
\begin{equation}
\delta w=\langle \delta \B,v_\B\rangle-\langle \delta \A,v_\A\rangle.
\end{equation}
By using the rotational invariance we may start with a Keplerian arc in a plane $\O xy$ with ends at equal ordinates $y_\A=y_\B$, and consider only variations in the same plane, having ends with the same property. As $\|\A\B\|$ is fixed, the pair $(\A,\B)$ is  translated and $\delta\A=\delta\B$. We get
\begin{equation}
\delta w=\langle \delta \A,v_\B-v_\A\rangle.
\end{equation}

\medskip\noindent{\bf Lemma 4.} Let $q\in\R^d$ be the position vector, $v=\dot q$ the velocity vector, $\varepsilon=q/\|q\|$ the radial unit vector. Consider two positions $\A$ and $\B$ on the same Keplerian orbit. The vectors $v_\B-v_\A$ and $\varepsilon_\A+\varepsilon_\B$ are linearly dependent.

\medskip\noindent{\bf Proof.} The proof in dimension $d=2$ or 3 with cross product notation will be more familiar to most readers (see remark 10). The eccentricity vector is
$$E=\varepsilon-v\times C, \quad\hbox{where }C=q\times v.$$
As we have $$0=\varepsilon_\A-\varepsilon_\B-(v_\A-v_\B)\times C,$$ the direction of $v_\A-v_\B$ is orthogonal to the direction of $\varepsilon_\A-\varepsilon_\B$, i.e., is the direction of $\varepsilon_\A+\varepsilon_\B$.\qed

\medskip\noindent{\bf End of proof of Proposition 1.} Now it is enough to prove that if  $\|\O\A\|+\|\O\B\|$ is fixed
\begin{equation}\label{old7}
\langle \delta \A,\varepsilon_\A+\varepsilon_\B\rangle=0.
\end{equation}
The vector $\delta \A$ is the variation of $\A$ relative to the fixed point $\O$. But we can also consider the relative variation of $\O$ in a translated frame where $\A$ and $\B$ are fixed. The relative infinitesimal variation of $\O$ is $\delta_{\rm rel}\O=-\delta \A$. But $\O$ is constrained to remain on a level set of the function $\O\mapsto \|\O\A\|+\|\O\B\|$, whose gradient is $-\varepsilon_\A-\varepsilon_\B$. This gradient is orthogonal to $\delta_{\rm rel}\O$, which proves (\ref{old7}).\qed

\medskip\noindent{\bf Proof of Theorem 1 (Hamilton).} Formula (\ref{3n}) suggests the following method to compute the elapsed time $\Delta t$ on any given nonrectilinear Keplerian arc. According to Lemma~3, the arc is not isolated among the arcs with same ends $\A$ and $\B$. Consider an infinitesimal variation $\delta$ among these arcs. Then $\Delta t=\delta w/\delta H$ according to (\ref{3n}), which means that $\Delta t$ is the derivative of $w$ with respect to $H$. Consider another variation of the arc, now with varying ends, such that $\|\A\B\|$, $\|\O\A\|+\|\O\B\|$ and $H$ are invariant. During this new variation,  $w$ is invariant according to Proposition 1, and thus $\Delta t$ is invariant as the derivative of an invariant with respect to another.\qed
\medskip

{\bf Some remarks about the proofs of (\ref{3n}).} In the variational calculus of Lagrange and Hamilton, $\delta$ is used together with $d$, which denotes the derivative along the trajectory. The commutation $d\delta=\delta d$ is stated and used (see \citealt{lagrange6}, seconde partie, \S IV, 3).  Concerning formula (\ref{3n}), which is his formula (A.), \citet{hamilton1} integrates by parts the middle term of his previous formula. He uses $d\delta=\delta d$, which would be incorrect if the trajectories were parametrized by the time $t$, since the interval of time is not constant during the variation. The trajectories are indeed reparametrized by $x$, $y$ or $z$. The editors of the {\it Mathematical papers} propose a similar proof after Hamilton's formula (Q.).

\medskip\noindent{\bf A proof of (\ref{3n}) by Jacobi.}  \citet{jacobi1}, \S II, proposed, as did \citet{hamilton2} just after his equation (28.), to use Hamilton's {\it principal function}, another ``action'' which is the integral of the Lagrangian:
\begin{equation}
\label{actionS}S=\int_{t_\A}^{t_\B} \bigl({1\over 2}\|v\|^2+U\bigr) dt,\quad\hbox{with }U={1\over r}.
\end{equation}
Here we specialize to the Kepler problem a general theory which is valid for the natural systems defined by a configuration space, a Riemannian metric and a potential $U$. The variation of $S$ among solutions starting all at the same time $t_\A$ and finishing all at the same time $t_\B$ is:
\begin{equation}
\delta S=\langle \delta \B,v_\B\rangle-\langle \delta \A,v_\A\rangle.
\end{equation}
This is simply the result of the integration by parts which gives the famous Lagrange equations of motion. We also have the more complete variation formula for solutions where $t_\A$ and $t_\B$ are allowed to vary:
\begin{equation}\label{4n}
\delta S=\langle \delta \B,v_\B\rangle-\langle \delta \A,v_\A\rangle-H(\delta t_\B-\delta t_\A).
\end{equation}
To get the new term,  consider the variation $\delta$ which consists in continuing the trajectory after the point $\B$. Due to (\ref{actionS}), $\delta S=(\|v_\B\|^2/2+U_\B)\delta t_\B$. We look for an expression of the form $\delta S=\langle \delta \B,v_\B\rangle-\langle \delta \A,v_\A\rangle+x(\delta t_\B-\delta t_\A)$.  Here $\delta\A=\delta t_\A=0$ but $\delta\B=(\delta t_\B) v_\B$, since the final position $\B$ is shifted along the trajectory. Consequently, $(\|v_\B\|^2/2+U_\B)\delta t_\B=\|v_\B\|^2\delta t_\B+x\delta t_\B$ and  $x=-\|v_\B\|^2/2+U_\B=-H$, which proves (\ref{4n}).

Along solutions, $H=\|v\|^2/2-U$ is constant and we have
$$S=\int_{t_\A}^{t_\B} (\|v\|^2-H) dt=w-H(t_\B-t_\A)$$
and 
$$
\delta S=\delta w-\delta H(t_\B-t_\A)-H(\delta t_\B-\delta t_\A).
$$
Comparing with (\ref{4n}), we get Hamilton's formula (\ref{3n}). \qed

\section{Elementary arguments}
\label{sect5}
The variational calculus offers considerable insights into Lambert's theorem. Hamilton's use of formula (\ref{3n}), when combined with the eccentricity vector (\ref{o8}), is a tool of remarkable flexibility which will be adapted to other problems in subsequent works. But in the present work we will set the variational ideas aside.

In the next section, we will show how the eccentricity vector can also be used to shorten elementary proofs of Lambert's theorem. Having already checked by a single differentiation that the quantities $\alpha$, $\beta$ and $C$ defined in (\ref{o8}) are constants of motion, we will not need differential calculus anymore. Our arguments will be as elementary as those of Lambert's original proof, and shorter. In this section, which is independent of Sects.\thinspace\ref{sect4} and \ref{sect6}, we show how the same elementary tools present themselves and give a  proof of Theorem 2.

\subsection{First attempt at proof of Theorem 2}
\label{sect5.1}
To describe the ``motion'' of the ends allowed by the hypothesis of Theorems 1 and 2, a good idea is to exchange what is fixed and what is moving. We represent $\O$ as moving relative to $\A$ and $\B$, as we did in the proof of Proposition 1 (see Fig.\thinspace6). The point $\O$ moves on an ellipse with foci $\A$ and $\B$. When $\A\neq\B$, it moves continuously until it reaches the line $\A\B$. There, a Keplerian arc with ends $\A$ and $\B$ should be rectilinear, since a nonrectilinear Keplerian orbit passes through the same point when it crosses a ray from $\O$ again.  {\it This argument will prove Theorem 2 as soon as we produce an arc with given energy $H$ for all the intermediate positions of $\O$.}

\vspace{0.5cm}
\centerline{\includegraphics[width=30mm]{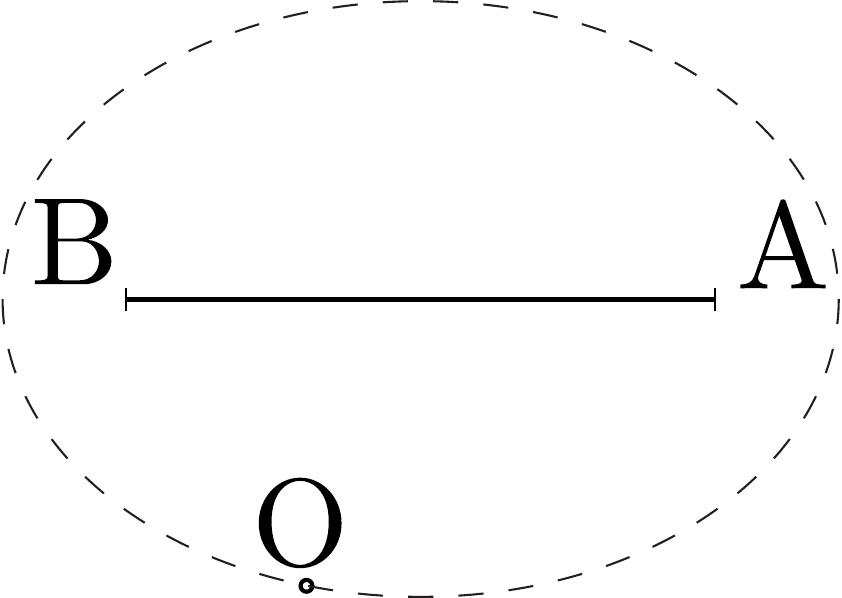}}
\centerline{\sl Fig.\thinspace6. How the focus $\O$ moves if  $\|\O\A\|+\|\O\B\|$ is constant.}
\vspace{0.5cm}

\medskip\noindent{\bf Proof for elliptic arcs.} Remark 3 suggests a proof of Theorem 2 in the elliptic case: in order to find the intermediate arcs, we find the second focus $\F$ by using Gauss's construction. The inequality $\|\O\A\|+\|\O\B\|+\|\A\B\|\leq 4a$, which guarantees that the two circles intersect, is satisfied during the continuous change since by hypothesis  $\|\O\A\|+\|\O\B\|$, $\|\A\B\|$ and $H=-(2a)^{-1}$ are constant. When $\O$ moves on Fig.\thinspace6, $\F$ moves smoothly and each successive position determines an intermediate orbit and an intermediate arc.\qed

\subsection{How the orbit moves}
\label{sect5.2}
While the point $\O$ describes an ellipse with foci $\A$ and $\B$ and major axis $\|\O\A\|+\|\O\B\|$, the point $\F$ describes an ellipse with same foci and major axis $4a-\|\O\A\|-\|\O\B\|$. If an orthogonal frame is chosen with origin the midpoint of $\A\B$, and $x$-axis along $\A\B$, then the abscissas of $\O$ and $\F$ remain in constant  proportion during the variation, and the same is true of their ordinates. These proportionalities are consequences of the identity $\|\O\B\|-\|\O\A\|=\|\F\A\|-\|\F\B\|$. In generic cases, either $\O$ and $\F$ are on the same side of $\A\B$, or they are on opposite sides (see Figs.\thinspace5, 7).

\vspace{0.5cm}
\centerline{\includegraphics[width=35mm]{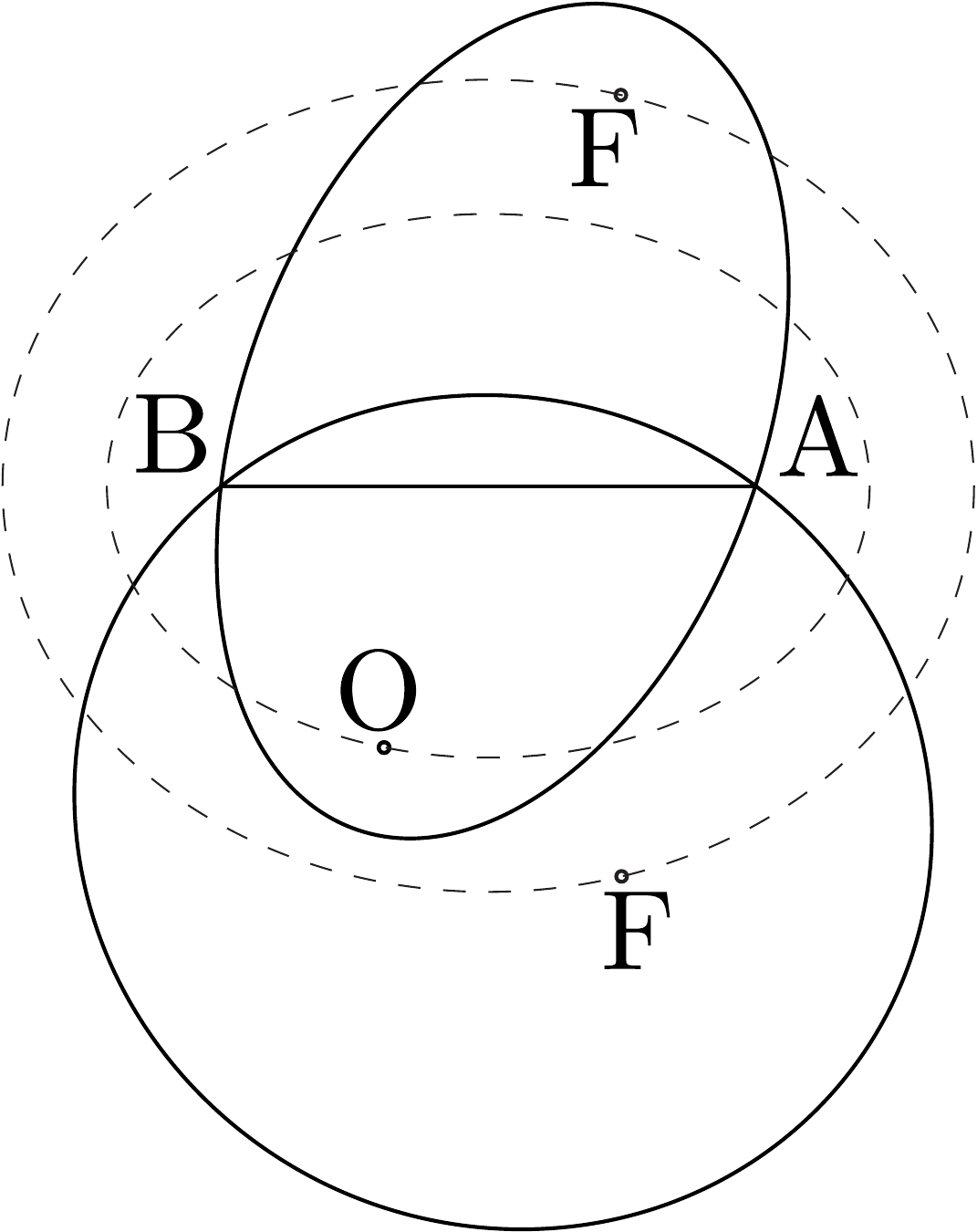}}
\centerline{\sl Fig.\thinspace7. How both foci $\O$ and $\F$ move if  $\|\O\A\|+\|\O\B\|$ and $H$ are constant.}

\subsection{Comment} 
\label{sect5.3}
The above elegant proof inspired by Gauss's construction is restricted to ellipses, while Theorems 1 and 2 do not distinguish the three kinds of conic section. Even if we can adapt the argument to hyperbolic orbits, and argue with a limiting process to get the parabolas, we should seriously wonder why we would need a proof for each kind of conic section. We will see that the continuous change of Theorem 2 may be described in the same way in the three cases.

We will argue with the eccentricity vector $E$ rather than the second focus. We recall that the consideration of $E$ was also a key to Proposition 1. In contrast with \citet{hamilton1} or with \citet{ioukovsky}, we were able to prove everything in a single case rather than three cases.

\subsection{Last remarks on variational calculus}
\label{sect5.4}
In the variational context of the previous section, a natural idea is to prove the existence of the intermediate arcs by minimizing an action integral.

Results related to the ``homogeneous calculus of variations'' are described in \citet{wintner}, \S 253, which uses \citet{todhunter}, \S 182.
They give the condition under which a Keplerian arc has a minimal Maupertuis action compared to the paths with same ends and same energy. They also determine the ``broken'' path which realizes the minimum when the condition is not satisfied. \citeauthor{wintner} recalls on p.\ 423 that the theory of the minimizing orbits is difficult.

We could also find the intermediate arcs through the nonhomogeneous calculus of variations. Having in mind Theorem 1 and remark 5, we are allowed to change the statement of Theorem 2: we ask for a path in ${\cal A}$ with constant $\|\A\B\|$, $\|\O\A\|+\|\O\B\|$ and $\Delta t$.  All along the path we need an arc with given $\A$ and $\B$ and with given $\Delta t$. This is the Lambert problem (see remark 7). We could get an arc by minimizing the action (\ref{actionS}).  However, the relevant techniques, for example such as presented in \citet{gordon}, are too sophisticated for an elementary result as Theorem 2. We would meet technical problems with the possibility of a collision orbit, and the smoothness of the variation would not be guaranteed.

\subsection{Rescaled Gauss's construction} 
\label{sect5.5}
Gauss's construction gives the position of the second focus $\F=(x_\F,y_\F)$ as an intersection of the circles with respective equations:
$$(x_\F-x_\A)^2+(y_\F-y_\A)^2=(2a-r_\A)^2,$$
$$(x_\F-x_\B)^2+(y_\F-y_\B)^2=(2a-r_\B)^2,$$
where $r_\A=\|\O\A\|=\sqrt{x_\A^2+y_\A^2}$, $r_\B=\|\O\B\|=\sqrt{x_\B^2+y_\B^2}$, and $a=(-2H)^{-1}$ is the semimajor axis. The equation of the line passing through both intersections is obtained by subtraction. As in the proof of Proposition 1, we set $y_\A=y_\B$. This equation is
$$x_\F={2a(r_\A-r_\B)\over x_\A-x_\B}.$$
If we set
$$\alpha={x_\F\over 2a},\quad \beta={ y_\F\over 2a},$$
then the position $(\alpha,\beta)$ of the ``rescaled second focus''  is an intersection of the line of equation
\begin{equation}\label{alpha}
\alpha={r_\A-r_\B\over x_\A-x_\B}
\end{equation}
with the circle of center $(x_\A/2a,y_\A/2a)$ and radius $1-r_\A/2a$ and with the circle of center $(x_\B/2a,y_\B/2a)$ and radius $1-r_\B/2a$. {\it We can now apply Gauss's construction to hyperbolas, and remarkably, to parabolas}, by simply replacing $-1/2a$ by $H$ in these formulas.
 
 In the parabolic case $H=0$, both circles coincide with the unit circle. The vertical line (\ref{alpha}) always determines two distinct points on the unit circle, since the triangular inequality gives $\alpha^2\leq1$, and $\alpha=\pm 1$ is the rectilinear case.
 
In the hyperbolic case the circles are tangent when $\pm(\|\O\A\|-2a)=\|\A\B\|\pm(\|\O\B\|-2a)$, where $a=-(2H)^{-1}$ is negative. Only the two cases with $\A$ and $\B$ on the same ray are possible. When $\A$ and $\B$ are not on the same ray, there are always two
distinct intersections: the number of intersections can change only where there is a tangency. 
 
Although the above discussion of the intersection of the circles is complete, we may ask for a direct discussion of the intersection of a circle with the vertical line (\ref{alpha}). Here are examples of formulas which are useful for this purpose.
$$\alpha={r_\A-r_\B\over x_\A-x_\B}={x_\A+x_\B\over r_\A+r_\B}={r_\A-r_\B+x_\A+x_\B\over  x_\A - x_\B +r_\A+r_\B},$$
\begin{equation}\label{alpha1}
\alpha-1={2(x_\B-r_\B)\over x_\A-x_\B+r_\A+r_\B}.
\end{equation}

\subsection{Proof of remark 2 and of Theorem 2}
\label{sect5.6} 
{\it The above ``rescaled second focus'' $(\alpha,\beta)$ is  the eccentricity vector $E$} defined in (\ref{o8}), since the norm of this vector is the eccentricity.
According to Gauss's rescaled construction, if we fix $H$ with the condition $H>H_{\min}$, where $H_{\min}$ is defined in (\ref{Hmin}), there are exactly two distinct values of the eccentricity vector $(\alpha,\beta)$ for which the Keplerian branch $r=\alpha x+\beta y+\gamma$ passes through two given points $\A$ and $\B$ which are not on the same ray. This corresponds to exactly two distinct Keplerian branches, since we have no choice for $\gamma$, which is determined by  $r_\A=\alpha x_\A+\beta y_\A+\gamma$. This proves the last claim in remark 2, which was proved in remark 3 in the case $H<0$.

As for the proof of Theorem 2, we should  find  a continuous change for which the arcs keep the same energy $H$ when $\O$ moves on the ellipse of Fig.\thinspace6. We simply follow an intersection of the moving circles and lines in Gauss's rescaled construction. The condition for intersection $H>H_{\min}$ is always satisfied, since $H_{\min}$ is constant when $\O$ moves.

There only remains to consider Theorem 2 in the particular case of a nonrectilinear arc with coinciding ends $\A=\B$. This is an elliptic orbit of eccentricity $e$ and semimajor axis $a$. We should push $e\to 1$ and keep $\A=\B$ on the orbit,  without changing $\|\O\A\|$ and $a$. This is always possible since as $e$ is increasing the interval $[a(1-e),a(1+e)]$ of the possible distances from the origin is increasing.\qed

\section{A short constructive proof of Theorem 1}
\label{sect6}
This theorem easily reduces to the bidimensional case $d=2$. We will prove it on $\R^2=\O xy$. This section is independent of Sects.\thinspace\ref{sect4} and \ref{sect5}, but uses again the unifocal equation presented in Sect.\thinspace\ref{sect3}.

In this proof the case $\A=\B$ requires a separate study. The velocity vectors $v_\A$  and $v_\B$ at times $t_\A$  and $t_\B$ have the same norm, which is given by the value of $H$.

If $v_\A=v_\B$,  the orbit is periodic, the point $q$ makes $k$ turns and comes back after a time $\Delta t=kT$ where the period $T$ is obtained through Kepler's third law: $T=2\pi a^{3/2}$. The elapsed time $\Delta t$ is constant along a path in ${\cal A}$ with constant $\|\A\B\|$ and constant $H=-(2a)^{-1}$. Theorem 1 is proved in this case. 

If $v_\A\neq v_\B$, then the only possibility is $v_\A=-v_\B\neq 0$. The orbit is rectilinear and goes back to the initial point after bouncing inward or culminating outward. In the space ${\cal A}$ of Keplerian arcs, this case defines, up to rotation, isolated arcs: $\Delta t$ is tautologically constant, which is Theorem 1 in this case.

We will assume $\A\neq \B$ in the rest of the proof. A curve in $\R^2$ is an irreducible conic section with a focus at $\O$ if and only if its equation is
\begin{equation}\label{o9}
x^2+y^2-(\alpha x+\beta y+\gamma)^2=0,\quad \hbox{for some  }(\alpha,\beta)\in\R^2,\;\gamma\in\;]0,+\infty[.\end{equation}
The eccentricity is $\sqrt{\alpha^2+\beta^2}$, the semiparameter is $\gamma$. Setting $r=\sqrt{x^2+y^2}\geq 0$, and excluding a branch in the hyperbolic case, we get  
$
r=\alpha x+\beta y+\gamma,
$
which is Eq.\thinspace(\ref{unifoc}).

We introduce an angular parameter $\phi$ and a number $f$. We write the identity $(x\cos\phi+f\sin\phi)^2+(x\sin\phi-f\cos\phi)^2=x^2+f^2$ in the form
\begin{equation}\label{o11}
(x\sin\phi-f\cos\phi)^2+y^2-f^2=x^2+y^2-(x\cos\phi+f\sin\phi)^2.
\end{equation}
\medskip\noindent{\bf Lemma 5.} For any $\phi\in\;]0,\pi[$, any $(\MU,\NU)\in\R\times\,]0,+\infty[$, the image of the Keplerian branch $\Sigma$ with equation $r=\MU y+\NU$ by the affine map
\begin{equation}\label{map2}
(x_1,y_1)\mapsto (x_2,y_2),\end{equation} where $$x_1=x_2\sin\phi-y_2\MU \cos\phi-\NU\cos\phi,\qquad y_1=y_2,$$
is a Keplerian branch around $\O$ with  equation $r=x\cos\phi+y\MU\sin\phi+\NU\sin\phi$.

\medskip\noindent{\bf Proof.}  In the equation $x_1^2+y_1^2-(\MU y_1+\NU)^2=0$, replace $(x_1,y_1)$ by its expression in $(x_2,y_2)$, then set $f=\MU y_2+\NU$ and use (\ref{o11}). \qed

\medskip
The principal axis of $\Sigma$ is vertical. If $\phi=\pi/2$ the affine map (\ref{map2}) is the identity. If not, (\ref{map2}) is not a transvection {\it stricto sensu} since $x_2$ has coefficient $\sin\phi\neq 1$. If we rescale by compounding with $(x_2,y_2)\mapsto (x_3,y_3)$, $x_2 = x_3/\sin\phi$, $y_2 = y_3/\sin\phi$, it is still not a transvection, now because of the coefficient $1/\sin\phi\neq 1$ of $y_3$. But it is  the transformation we need.

\medskip\noindent{\bf Lemma 6.} For any $\phi\in\;]0,\pi[$, any $(\MU,\NU)\in\R\times\,]0,+\infty[$, the image of the Keplerian branch $\Sigma$ with equation $r=\MU y+\NU$ by the affine map
\begin{equation}\label{o12}
(x_1,y_1)\mapsto (x_3,y_3),
\end{equation}
where
$$
x_1=x_3-\MU {\cos\phi\over \sin\phi}\,y_3-\NU\cos\phi,\qquad y_1={1\over \sin\phi}\,y_3,
$$
is a Keplerian branch around $\O$ with equation $r=x\cos\phi+y\MU\sin\phi+\NU\sin^2\phi$.

\medskip\noindent{\bf Lemma 7.} To any Keplerian branch $\Omega$ with equation $r=\alpha x+\beta y+\gamma$,  $\gamma>0$, which possesses a horizontal chord, is uniquely associated a triple $(\phi, \MU,\NU)\in\;]0,\pi[\times\R\times]0,+\infty[$ such that $\Omega$ is the image by the affine map (\ref{o12}) of the vertical branch $\Sigma$ with equation $r=\MU y+\NU$.

\medskip\noindent{\bf Proof.} There are two points $\A$ and $\B$ of same ordinate $y_\A=y_\B$ on the branch $\Omega$. Let $r_\A = \|\O\A\|$, $r_\B = \|\O\B\|$. Then $r_\A = \alpha x_\A + \beta y_\A + \gamma$,  $r_\B = \alpha x_\B + \beta y_\B + \gamma$  give $r_\A - r_\B=\alpha (x_\A - x_\B)$. The triangular inequality gives
 $\alpha^2<1$. We set $\phi=\arccos\alpha$, $\MU=\beta/\sin\phi$, $\NU=\gamma/\sin^2\phi$.\qed

\medskip\noindent{\bf Lemma 8.} The affine map (\ref{o12}) sends any horizontal chord $\A\B$ of the branch $\Sigma$ onto a horizontal chord of same length $\|\A\B\|$ and same $\|\O\A\|+\|\O\B\|$.

\medskip\noindent{\bf Proof.} The unit coefficient of $x_3$ in formula (\ref{o12}) gives the invariance of the horizontal length $\|\A\B\|$. Now as in the previous proof, $r_\A = \alpha x_\A + \beta y_\A + \gamma$ and $r_\B = \alpha x_\B + \beta y_\B + \gamma$ give
by subtraction $\alpha=(r_\A-r_\B)/(x_\A-x_\B)$ which is also $(x_\A+x_\B)/(r_\A+r_\B)$. By addition $r_\A+r_\B=\alpha^2(r_\A+r_\B)+2\beta y+2\gamma$. Replacing $\alpha=\cos\phi$, $\beta=M\sin\phi$, $\gamma=N\sin^2\phi$, this is
$$r_\A+r_\B={2\MU y\over \sin\phi}+2\NU.$$ Here $y$ is $y_3$. As $y_3=y_1\sin \phi$, $r_\A+r_\B=2My_1+2N$ does not vary with $\phi$.\qed
\medskip

After these geometrical lemmas we pass to dynamics. We consider with slight abuse of terminology that the affine map (\ref{o12}) sends a Keplerian orbit onto a Keplerian orbit and a Keplerian arc onto a Keplerian arc. The map indeed sends an unparametrized arc onto an unparametrized arc.  The running direction is induced by the affine map, but the time parametrization is not. Instead, we consider on the image branch a time parametrization compatible with Newton's system (\ref{n1}).

\medskip\noindent{\bf Lemma 9.}  Under the affine map (\ref{o12}), the areal constant $C$ is shrunk as the areas, i.e., is multiplied by $\sin\phi$, while the energy $H$ is unchanged.

\medskip\noindent{\bf Proof.} The affine map sends a conic section with semiparameter $\NU$ to a conic section with semiparameter $\NU\sin^2\phi$. As $C=\pm\sqrt{\gamma}$, the angular momentum $C$ is multiplied by $\sin\phi$. The Jacobian of the map is also $\sin\phi$, as the ordinates are multiplied by  $\sin\phi$ and the horizontal distances are preserved. The left-hand side of  (\ref{H}) was $\MU^2-1$ and becomes $\cos^2\phi+\MU^2\sin^2\phi-1=(\MU^2-1)\sin^2\phi$.\qed

\medskip\noindent{\bf Lemma 10.}  The elapsed time $\Delta t=t_\B-t_\A$ on a Keplerian arc with endpoints $\A$ and $\B$ having the same ordinate $y_\A=y_\B$ is unchanged by the affine map (\ref{o12}). 

\medskip\noindent{\bf Proof.} The period $T=2\pi (-2H)^{-3/2}$ of an elliptic orbit is unchanged, since the energy $H$ is unchanged. This reduces the study to arcs making less than one turn around $\O$. The elapsed time is twice the area swept out divided by the areal constant. According to Lemma 9, it is enough to prove that the area swept out is shrunk as are all the areas, i.e., is multiplied by $\sin\phi$. But this area is made of a segment of conic section ${\cal S}$ and of the triangle ${\cal T}=\O\A\B$ (see Fig.\thinspace8). The segment ${\cal S}$ is sent onto the corresponding segment, whose area is shrunk as is any area. The triangle ${\cal T}$ is sent onto a triangle formed by the image of the chord $\A\B$ and the image of $\O$, which is on the $x$-axis. Its area is shrunk as is any area. This is not the triangle we need in order to compute the area swept out from $\O$. But the triangle we need, with vertex at $\O$, has the same area. \qed
\medskip

We summarize the results of this section in a proposition.

\medskip\noindent{\bf Proposition 2.}  Starting from any given Keplerian arc $\Gamma$, we can arrive at some rectilinear arc by a continuous change which keeps constant the four quantities $\|\A\B\|$, $\|\O\A\|+\|\O\B\|$, $H$ and $\Delta t$. This change may be described as follows. In a frame $\O xy$, the chord $\A\B$ remains horizontal. Let $\G$ be its midpoint, let $\phi$ be the parameter of the change, $\rho$, $\sigma$ be constants such that $\O\G=(\rho\cos\phi,-\sigma\sin\phi)$. The angle $\phi$ varies from a $\phi_\Gamma\in\,] 0,\pi[$   to $\phi=0$. The equation of the  branch carrying the arc is $r=x\cos\phi+y\MU\sin\phi+\NU\sin^2\phi$. The constants $\rho$, $\sigma$, $\phi_\Gamma$, $\MU$, $\NU$ are smooth functions of $\Gamma$, satisfying the compatibility condition $\rho+\MU \sigma=\NU$.
\medskip
\vspace{0.5cm}
\centerline{\includegraphics[width=65mm]{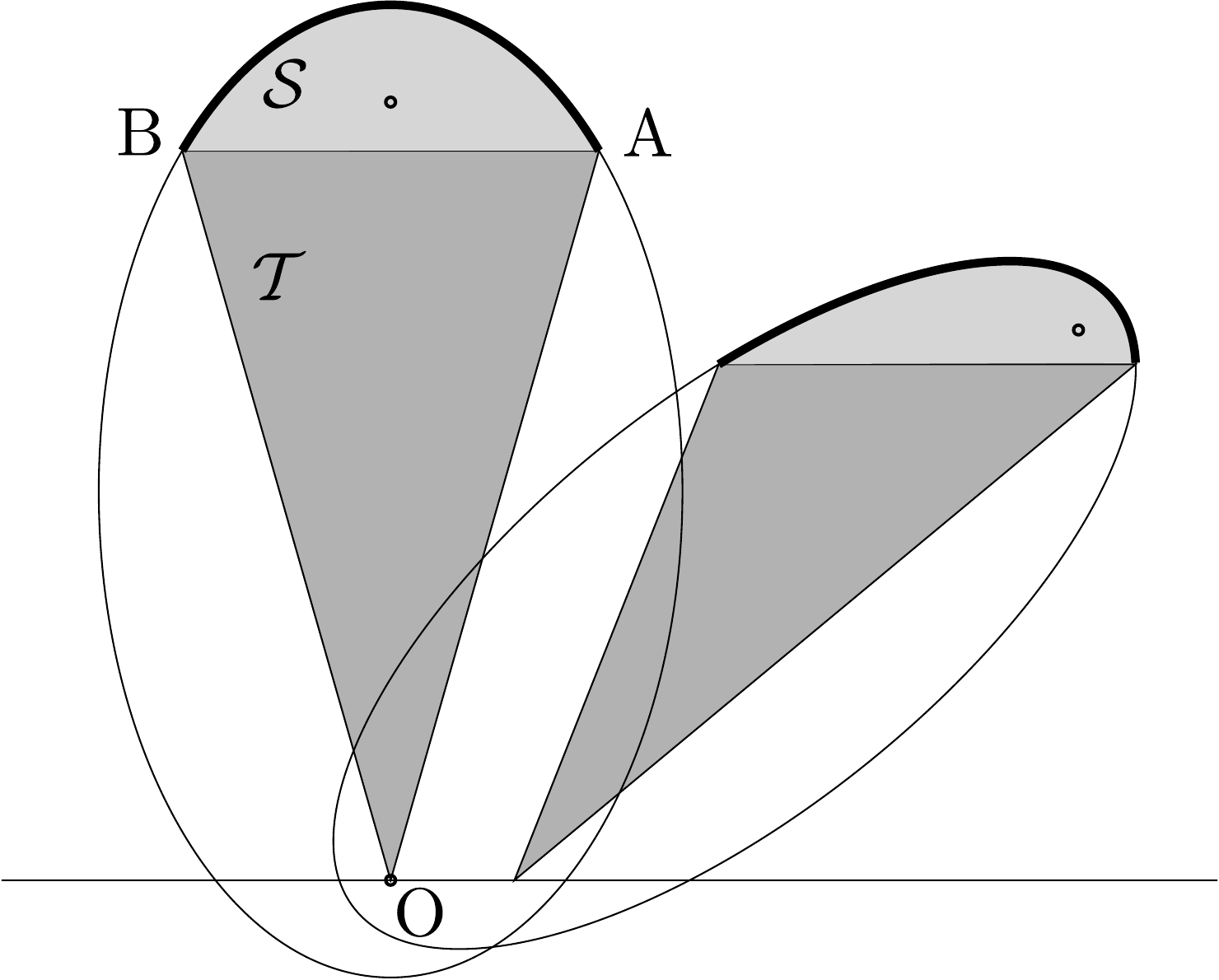}}
\centerline{\sl Fig.\thinspace8. Image of the swept out domain by the affine map (\ref{o12}).}
\vspace{0.5cm}

This statement is stronger than Theorem 2 since it also claims that $\Delta t$ is constant. It seems weaker than Theorem 1 since the latter is about an arbitrary continuous change with constant $(\|\A\B\|, \|\O\A\|+\|\O\B\|,H)$. But Theorem 1 is an easy corollary. If we start from such an arbitrary change, we can {\it project} it by sending each arc to the rectilinear arc which ends the continuous change described in Proposition 2. The projected change is continuous. It is constant, since  $(\|\A\B\|, \|\O\A\|+\|\O\B\|,H)$ determines a discrete choice of rectilinear arcs. Thus $\Delta t$ is constant. {\it Theorem 1 is proved again.}

\section{Noncrossing at rectilinear arcs}
\label{sect7}
The following proposition, which strengthens Theorem 2, is not completely established by the arguments in Sects.\thinspace\ref{sect5} and \ref{sect6}.

\medskip\noindent{\bf Proposition 3.} In the space ${\cal A}'$ of Keplerian arcs $\Gamma$ in the plane $\O xy$, having distinct ends $\A$ and $\B$ with same ordinate $y_\A=y_\B$, the nonempty connected components of the level sets of the map ${\cal A}'\to \R^3$, $\Gamma\mapsto (\|\A\B\|, \|\O\A\|+\|\O\B\|,H)$ are topologically circles. Each such circle contains two rectilinear arcs.
 
\medskip\noindent{\bf Remark 12.} Figure 6 suggests this proposition. Figure 7 seems to confirm it for $H<0$. The topological circles in the statement would correspond to an $\O$ describing a complete ellipse, while $\F$ describes its own ellipse, being the image of $\O$ by some affinity. Let us restate Proposition 3 according to this remark.

\medskip\noindent{\bf Definition 4.} A {\it Lambert cycle} of planar Keplerian arcs consists of 
\begin{itemize}[noitemsep,topsep=2pt,leftmargin=24pt]
\item[(i)]an arc $\Gamma$ carried by a vertical branch $\Sigma$ with equation $r=\MU y+\NU$ with $\NU>0$, with two distinct endpoints at the same ordinate,
\item[(ii)]all the images of $\Gamma$ by the affine maps (\ref{o12}), for all $\phi\in\;]0,\pi[$,
\item[(iii)]all the reflected arcs with respect to the horizontal axis, which indeed correspond to the $\phi\in\;]-\pi,0[$,
\item[(iv)]both limiting rectilinear Keplerian arcs as $\phi\to 0$ and as $\phi\to\pi$.
\end{itemize}

\medskip\noindent{\bf Proposition 4.} The connected components  described in Proposition 3 are the Lambert cycles.

\medskip\noindent{\bf Remark 13.} Remark 12 is based on the description in Sect.\thinspace\ref{sect5.2}. Definition 4 uses the affine maps of Sect.\thinspace\ref{sect6}. Let us show the compatibility of both descriptions. Call $\G$ the midpoint of $\A\B$. The vector $\O\G=(\rho\cos\phi,-\sigma\sin\phi)$, where $2\rho=\|\O\A\|+\|\O\B\|$, $\sigma^2=\rho^2-c^2$, $2c=\|\A\B\|$. The vector $\G\F=\bigl((2a-\rho)\cos\phi,\sigma_\F\sin\phi\bigr)$, where $\sigma_\F^2=(2a-\rho)^2-c^2$. Proposition 2 gives $\O\F=2a(\cos\phi,\MU\sin\phi)$. We check that $\O\G+\G\F=\O\F$ by checking that $2a\MU=\sigma_\F-\sigma$ is the relation between the eccentricity vector and the ordinates of the foci when $\phi=\pi/2$.

\medskip\noindent{\bf Proof.} The end of this section constitutes a proof of Propositions 3 and 4. If a connected component of a level set of $(\|\A\B\|, \|\O\A\|+\|\O\B\|,H)$ contains a nonrectilinear arc, it contains the Lambert cycle which passes through this arc. This cycle is computed by using Lemma 7 and Definition 4. Gauss's rescaled construction shows the local uniqueness of the arc in the level set, above each position of $\O$ in Fig.\thinspace6, until $\O$ reaches the line $\A\B$. So, if two Lambert cycles are in the same level set, they can only connect at a rectilinear arc.

May a rectilinear arc belong to several Lambert cycles? We will answer negatively by  counting the rectilinear arcs and the Lambert cycles. In Figs.\thinspace9 and 10, a rectilinear orbit with two marked points $\A$ and $\B$ is approached by two distinct paths of marked orbits with same $(\|\A\B\|, \|\O\A\|+\|\O\B\|,H)$. In Fig.\thinspace9, the focus $\F$ is above the chord $\A\B$, while in Fig.\thinspace10, $\F$ is below $\A\B$. We can also distinguish these two choices in Fig.\thinspace7. They are the only choices, as shown by Gauss's construction.

Consider the rectilinear arcs from $\A$ to $\B$ for which $\Delta t$ is shorter than a period. We can see four such arcs with same energy, and the efficient way to describe them is in terms of these general cases:

\medskip\noindent{\bf Definition 5.} A Keplerian arc around $\O$ making less than one turn is said to be
\begin{itemize}[noitemsep,topsep=2pt,leftmargin=18pt]
\item[--]indirect, or $I_\O$, if its convex hull contains $\O$; direct, or $D_\O$, if its convex hull does not contain $\O$;
\item[--]indirect with respect to the second focus $\F$, or $I_\F$, if its convex hull contains $\F$; direct with respect to $\F$, or $D_\F$, if its convex hull does not contain $\F$.
\end{itemize}
\medskip
\vspace{0.5cm}
\centerline{\includegraphics[width=56mm]{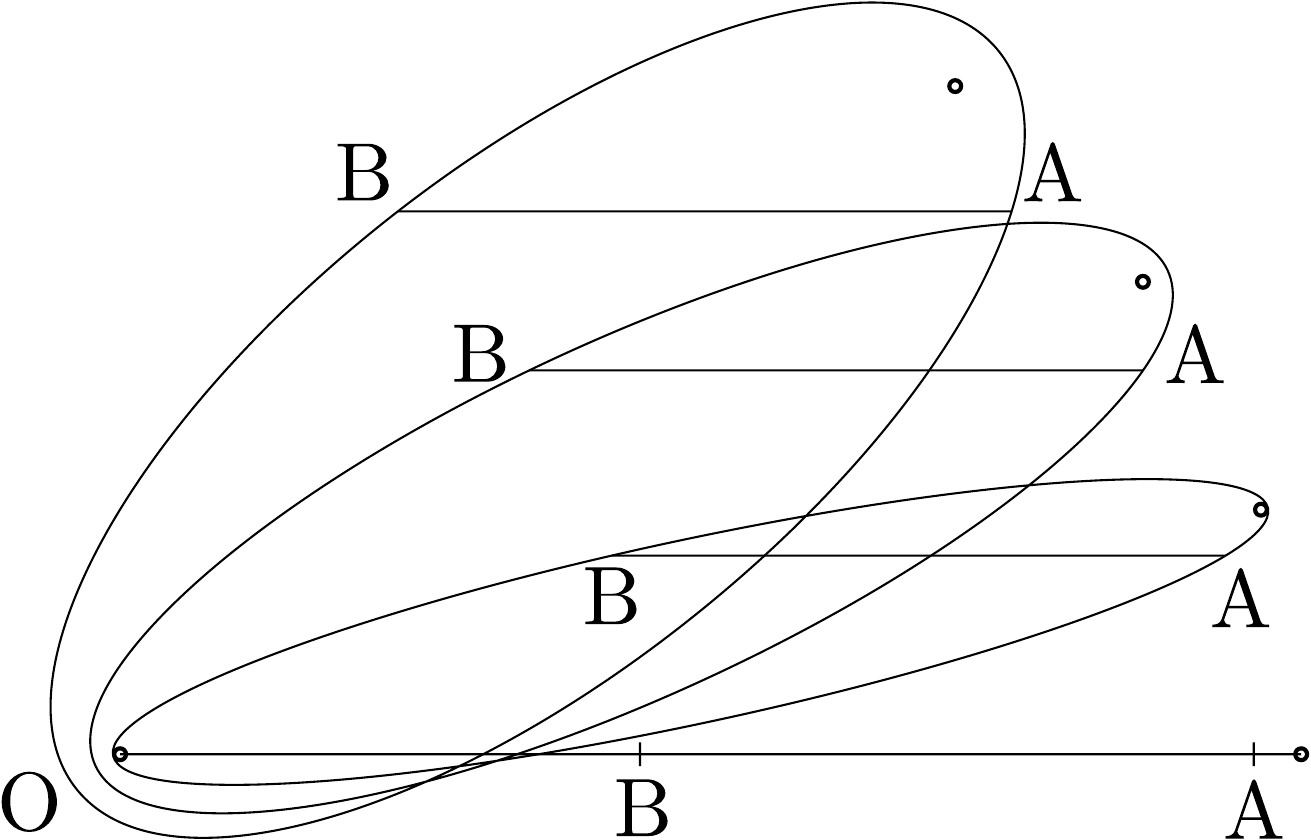}}
\centerline{\sl Fig.\thinspace9. Ellipses with a chord, with same $\|\A\B\|$, $\|\O\A\|+\|\O\B\|$ and $H$.}
\vspace{0.5cm}

In Figs.\thinspace9 and 10, the culmination point of a rectilinear elliptic arc is the limiting second focus $\F$. When a focus belongs to the convex hull of an arc,
this property is preserved in the limiting rectilinear arc. In other words, we have:

\medskip\noindent{\bf Proposition 5.} If a Keplerian arc in a Lambert cycle is $D_\O$ (respectively $I_\O$, $D_\F$, $I_\F$)  all the Keplerian arcs of the cycle are $D_\O$ (respectively $I_\O$, $D_\F$, $I_\F$).
\medskip

\vspace{0.5cm}
\centerline{\includegraphics[width=60mm]{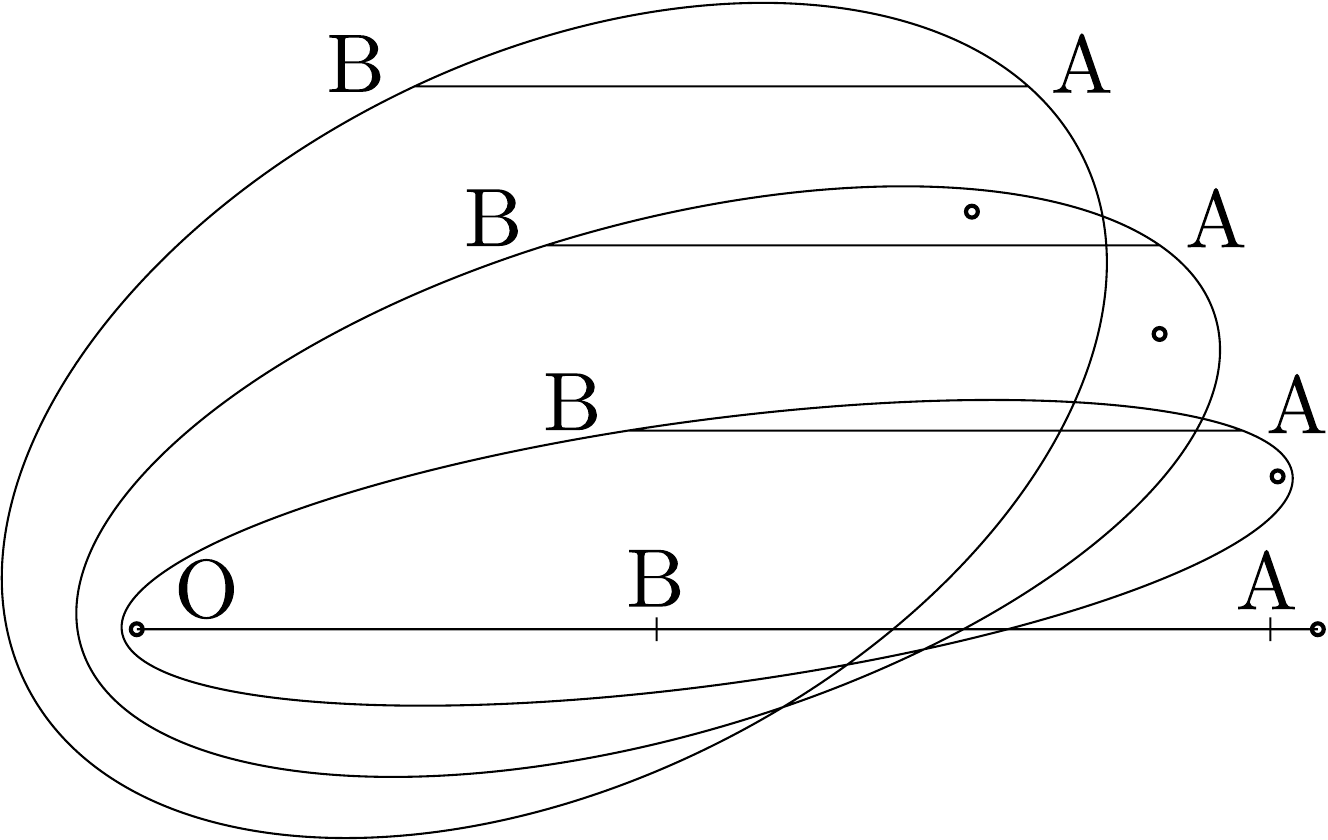}}
\centerline{\sl Fig.\thinspace10. Another family, with same chords and same energy as in Fig.\thinspace9.}
\vspace{0.5cm}

On a given rectilinear arc the body collides or does not,  and independently, culminates or does not. To say it in another way, on a given arc, the final velocity $v_\B$ belongs to $]0,+\infty]$ or to $]-\infty,0[$, and independently,  $v_\A\geq 0$ or $v_\A<0$. Together with Definition 5, we have three equivalent criteria distinguishing  four rectilinear arcs of same negative energy, the shortest being $D_\O D_\F$, the longest $I_\O I_\F$, the other two  $I_\O D_\F$ and $D_\O I_\F$.

The $I_\O D_\F$ and the $D_\O I_\F$ arcs are approached by nonrectilinear arcs in Fig.\thinspace9, but not in Fig.\thinspace10, while the $D_\O D_\F$ and the $I_\O I_\F$ arcs are approached in Fig.\thinspace10, but not in Fig.\thinspace9. There is a unique such approach of each of the four rectilinear arcs. In other words, {\it each of the four rectilinear arcs belongs to a unique Lambert cycle}. If the rectilinear arc is $I_\O D_\F$ or $D_\O I_\F$, the approach is by orbits with $\O$ and $\F$ separated by the chord. If it is
$D_\O D_\F$ or $I_\O I_\F$, $\O$ and $\F$ are on the same side  of the  chord.

The parabolic and hyperbolic arcs are always $D_\F$. There are only two types of arcs, $D_\O$ and  $I_\O$. During an approach of a rectilinear arc with $H\geq 0$, Gauss's rescaled construction gives two choices of eccentricity vector for each configuration $\O \A\B$, each giving in turn an orbit and an arc. A choice approaches the $I_\O$ rectilinear arc, the other the $D_\O$ rectilinear arc. Again, {\it each of both rectilinear arcs belongs to a unique Lambert cycle.}

We leave to the reader the interesting study of the limiting cases with a focus on the boundary of the convex hull of the arc. We should however make clear that a Lambert cycle is topologically a circle. This is obvious if $\O$ or $\F$ describes a nondegenerate ellipse in Fig.\thinspace7, or if the eccentricity vector $E$ describes a nondegenerate ellipse in the plane $\O xy$. Only in one case $\O$, $\F$ and $E$ all describe flat ellipses. This is when the foci $\O$ and $\F$ describe in Fig.\thinspace7 the segment $\A\B$, remaining opposite to each other. The full cycle has twice the same position of $(\O,\F)$, but the arcs are indeed distinct, being once the upper arc and once the lower arc, according to (\ref{o12}). So, even in this case, the family of arcs is topologically a circle in the space of arcs.

\section{Geometrical analogs}
\label{sect8}
Some geometrical statements are closely related to Lambert's theorem. Even if they can be expressed in many simple ways, they do not appear to be well known. We will give three propositions. We begin with a lemma published in \citet{terquem} as Theorem V.

\medskip\noindent{\bf Lemma 11.} In an ellipse
\begin{itemize}[noitemsep,topsep=2pt,leftmargin=18pt]
\item[--] an arbitrary chord passing through a focus,
\item[--]  the parallel chord passing through the center and
\item[--] the major axis
\end{itemize}
have their three lengths in geometric progression.

\medskip\noindent{\bf Proof.} We take the direction of the chord as the $x$-axis and a focus as the origin. We compute the horizontal semichord at ordinate $y$ as $\sqrt{\delta}/(1-\alpha^2)$ where $\delta=\alpha^2(\beta y+\gamma)^2-(1-\alpha^2)(y^2-(\beta y+\gamma)^2)=(\alpha^2-1)y^2+(\beta y+\gamma)^2$ is the reduced discriminant of Eq.\thinspace(\ref{o9}) seen as a trinomial in $x$. At $y=0$ the semichord is $\gamma/(1-\alpha^2)$. At $y=\beta a$, where $a$ is the semimajor axis, which satisfies $\gamma=a(1-\alpha^2-\beta^2)$,  $\delta=(\alpha^2-1)a^2\beta^2+a^2(1-\alpha^2)^2=a(1-\alpha^2)\gamma$. Consequently $\sqrt{\delta}/(1-\alpha^2)$ is the geometric mean of $a$ and $\gamma/(1-\alpha^2)$.\qed

\medskip\noindent{\bf Proposition 6.} Consider in a Euclidean plane an ellipse and a chord. Apply an affine map. Any two of these three properties imply the remaining one:
\begin{itemize}[noitemsep,topsep=2pt,leftmargin=18pt]
\item[--] a parallel chord passing through a focus is sent onto a chord passing through a focus,
\item[--] the length of the given chord is preserved,
\item[--] the length of the major axis is preserved.
\end{itemize}

\medskip\noindent{\bf Proof.} Call the three lengths in Lemma 11, corresponding to the direction of the given chord, $f$, $g$, $h$ before applying the map and $f'$, $g'$, $h'$ after. An affine map sends all the parallel chords to parallel chords, multiplying their length by a common factor $\lambda$. As the center of the ellipse is sent to the center of the image, we have $g'=\lambda g$. Observe now that in the family of parallel chords the length starts from zero, increases until it reaches a maximum and then decreases to zero.  Thus a ``parallel chord passing through a focus'' is also a ``parallel chord of same length as a parallel chord passing through a focus'': the first condition in the statement is $\lambda f=f'$. The second is $g=g'$, i.e., $\lambda=1$, the third is $h'=h$.
An easy analysis shows that if two conditions are satisfied, the geometric progression implies the remaining one. \qed
\medskip

Lemma 6 provides affine maps satisfying the three properties in Proposition 6, if the term $-\NU\cos\phi$ is removed from expression (\ref{o12}). But Lemma 6 works as well for parabolas and hyperbolas. Let us extend Proposition 6 accordingly.

\medskip\noindent{\bf Proposition 7.} The image of a conic section with semiparameter $\gamma>0$, with a focus on a straight line $D$, by an affine map with Jacobian determinant $J$, which fixes all the points of $D$, has semiparameter $J^2\gamma$ if and only if it has a focus on $D$.

\medskip\noindent{\bf Remark 14.} In the hypothesis and in the conclusion we should consider that a parabola has a focus at infinity, which is a point on the line at infinity. This focus is on $D$ if and only if the axis of the parabola is parallel to $D$. The proof of this proposition is a case-by-case study, which we leave to the reader.
\medskip

In the case of an ellipse or a hyperbola, the image has semiparameter $J^2\gamma$ if and only if the semimajor axis $a$ is preserved. To prove this we may use the expression $\pi a^{3/2}\sqrt{\gamma}$ of the area of an ellipse, and the expression $|a|^{3/2}\sqrt{\gamma}$ of the area of a triangle delimited by the two asymptotes of a hyperbola and a tangent.

\medskip\noindent{\bf Proposition 8.} Consider an ellipse drawn in an affine plane. Consider two Euclidean forms making this plane Euclidean in two different ways, each defining a pair of foci, each defining a major axis of the ellipse. Both major axes have equal length if and only if both Euclidean forms induce equal units of length on a chord passing through two foci, one of each pair.

\medskip\noindent{\bf Proof.} Consider that the first Euclidean form defines the Euclidean structure of the plane, and that the second is the pull-back of the first by an affine map.  By the well-known theory of the Gram matrix, there exists an affine map with such a pull-back. Apply Proposition 6 to the ellipse and a chord passing through two foci, one of each pair. The first hypothesis of Proposition 6 is satisfied. The second and the third are then equivalent.\qed

\section{A plethora of demonstrations}
\label{sect9}
Here is a timeline of Lambert's theorem.

{\bf 1687.} Newton considers the problem of the determination of the orbit of a comet from three observations, in Proposition XLI, Book III, of his {\it Principia}. \citet{lagrange3, lagrange6} will show later how two of Newton's lemmas give a proof of Lambert's theorem in the parabolic case (see \citealt{kriloff1, kriloff2}).

{\bf 1743.} \citeauthor{euler1} considers the same problem as Newton, about comets on parabolic orbits, and concludes \S XIII by the formula:
\begin{equation}\label{Euler}
6(t_\B-t_\A)=\bigl(\|\O\A\|+\|\O\B\|+\|\A\B\|\bigr)^{3/2}-\bigl(\|\O\A\|+\|\O\B\|-\|\A\B\|\bigr)^{3/2}.
\end{equation}
The time $t$ from the collision to the position $x$ in a rectilinear Keplerian motion with zero energy satisfies $6t=(2x)^{3/2}$. Euler's formula is this expression together with the reduction to the rectilinear case proposed by Theorems 1 and 2. Euler's choice of sign happens to correspond to a direct arc. His proof, which we call ${\cal P}_1$, is based on a simplification that appears when dividing the area swept out by the square root of the semiparameter. In \S XIV, Euler gives another proof, of the same nature, of the same formula. In \S XV, he considers the ``more difficult'' elliptic case and gets a formula in terms of the eccentricity and the three distances which is not as elegant. 

{\bf 1744.} At the opportunity of observing another comet, \citeauthor{euler2} reconsiders the determination of nearly parabolic orbits in a book.  He presents the computations differently and does not mention Eq.\thinspace(\ref{Euler}). He compares his new numerical results with what he got in 1743.

{\bf 1761.} In a letter in February (see \citealt{bopp}), Lambert announces to Euler his discovery of formula (\ref{Euler}).

{\bf 1761.} \citeauthor{lambert} publishes his fundamental book where he presents formula (\ref{Euler}) in \S 63, giving a proof of style ${\cal P}_1$, and later the elliptic case of our Theorem 1. His main step is the construction, from a general ellipse with a chord, of another ellipse with a chord perpendicular to the principal axis. This corresponds to our Lemma 7, the second ellipse being our vertical branch $\Sigma$.  Our Fig.\thinspace8 should be compared to Lambert's figure~21, our Lemma 9, to his \S 173, our Lemma 8, to his \S 177, our Lemma 10, to his \S 178.  Lambert expresses, through the rectilinear motion, the elliptic $\Delta t$ as an integral and as a series in \S 210, which he uses in \S 211 to obtain again formula (\ref{Euler}) as a limiting case. He mentions the hyperbolic case in \S 213, but only about the rectilinear motion. We will number this proof ${\cal P}_2$. Lambert refers to Euler's book of 1744 in his introduction, but not to Euler's article of 1743.

{\bf 1761.} In March Lambert sends his book to Euler who answers ``Votre theoreme pour exprimer l'aire d'un secteur parabolique est excellent, j'en puis bien voir la verit{\'e}, mais par de tels detours, que je ni serois jamais arriv{\'e}, si je ne l'avois su d'avance\thinspace; je attend donc avec impatience de voir l'analyse qui y a conduit sans detour\footnote{Your theorem for expressing the area of a parabolic sector is excellent, I can see the truth of it, but by such detours, that I could never have arrived to it had I not known it in advance; I therefore wait impatiently to see the analysis leading to it without detours.}'' in a first letter, ``la belle demonstration de l'aire du secteur parabolique, dont Vous m'avies communiqu{\'e} l'expression m'a caus{\'e} un tr{\`e}s sensible plaisir\thinspace; mais je fus bien plus surpris d'en voir l'application aux secteurs elliptiques [...] je reconnois aisement que les methodes, que j'avois propos{\'e}es autrefois, peuvent etre tr{\`e}s considerablement perfectionn{\'e}es\footnote{the beautiful proof of the area of a parabolic sector, the expression for which You communicated to me gave me great pleasure; but I was even more surprised to see its application to elliptic sectors [...] I easily recognize that the methods I proposed earlier may be improved considerably.}''  in a second, ``Vos remarques sur la reduction du mouvement curviligne des corps celestes {\`a} la chute rectiligne sont tr{\`e}s sublimes, et nous decouvrent en effet des prome[ss]es qui sans cette reduction paroissent tout {\`a} fait indechiffrables\footnote{Your remarks on the reduction of curvilinear motion of celestial bodies to rectilinear fall are very sublime, and we discover in fact promises which, without that reduction, appear indecipherable.}'' in a third.

{\bf 1773.} \citeauthor{lagrange2} deduces Lambert's theorem while discussing Euler's two fixed centers problem and analyzing the limiting case where one of the centers has zero mass and is on the orbit. See \S XI. We call this proof ${\cal P}_3$. See 1780, 1815, Jacobi 1866. We will briefly discuss in Sect.\thinspace\ref{sect10.1} a related work by \citeauthor{lagrange1} published in the same volume of {\it Miscellanea Taurinensia}.

{\bf 1780.} In memory of his friend and colleague who died in 1777, Lagrange publishes a series of memoirs. In the \citetalias{lagrange3}, he writes ``C'est ce que M.\ Lambert a fait depuis dans son beau Trait{\'e} {\it De orbitis Cometarum}, o{\`u} il est parvenu {\`a} un des Th{\'e}or{\`e}mes les plus {\'e}l{\'e}gants et les plus utiles qui aient {\'e}t{\'e} trouv{\'e}s jusqu'ici sur ce sujet, et qui a en m{\^e}me temps l'avantage de s'appliquer aussi aux orbites elliptiques\footnote{This is what Mr Lambert has since done in his beautiful Treatise {\it De orbitis Cometarum}, where he arrived at one of the most elegant and useful theorems produced on the subject up to now, with the additional advantage of applying also to elliptical orbits.}'' and ``Th{\'e}or{\`e}me qui, par sa simplicit{\'e} et par sa g{\'e}n{\'e}ralit{\'e}, doit {\^e}tre regard{\'e} comme une des plus ing{\'e}nieuses d{\'e}couvertes qui aient {\'e}t{\'e} faites dans la Th{\'e}orie du syst{\`e}me du monde\footnote{A Theorem which, by its simplicity and generality, must be regarded as one of the most ingenious discoveries made in the Theory of the system of the world.}''. He also analyses Euler's book of 1744, and several published consequences of Lambert's theorem.

 {\bf 1780.} \citeauthor{lagrange4} presents three other proofs of Lambert's theorem, introducing them in \S 1 by ``mais ce th{\'e}or{\`e}me m{\'e}rite particuli{\`e}rement l'attention des G{\'e}om{\`e}tres par lui-m{\^e}me, et parce qu'il parait difficile d'y parvenir par le calcul\thinspace; en sorte qu'on pourrait le mettre dans le petit nombre de ceux pour lesquels l'Analyse g{\'e}om{\'e}trique semble avoir de l'avantage sur l'Analyse alg{\'e}brique\footnote{but this theorem particularly deserves the attention of Geometers by itself, and because it seems difficult to achieve by calculation; so that one may place it among the small number of those for which Geometric analysis seems to have an advantage over Algebraic analysis.}.''  He rejects his proof ${\cal P}_3$ as too indirect and complicated, but proposes a similar proof ${\cal P}_6$ which does not refer explicitly to the two fixed centers problem (see \S 14). The first proof in \citet{lagrange4}, which we call ${\cal P}_4$, uses the eccentric anomaly. Note that the difference $u_\B-u_\A$ of the final and initial eccentric anomalies is obviously an invariant of our map (\ref{o12}), and that $(u_\A+u_\B)/2$ is the eccentric anomaly of the highest or lowest point of the ellipse. The second proof, which we call ${\cal P}_5$, is concluded in \S 7.  It starts with the expression of the elliptic $\Delta t$ by a quadrature of a function of the distance $r$, and then uses general methods rather than formulas for the Keplerian motion. In the proof ${\cal P}_6$, the concluding  identity of  ${\cal P}_5$ is presented as a particular case of more general identities (see our Sect.\thinspace\ref{sect10.1}).

{\bf 1784.} \citeauthor{lexell}, in a volume announcing the death of his master Euler, discusses the proofs by Lambert and Lagrange, extends them to the case of hyperbolic motions, and discusses reality conditions in Lagrange's identities. He also proposes some reciprocal statements.

{\bf 1797.} \citeauthor{olbers} publishes a method of orbit determination in a book with many references, including to formula (\ref{Euler}), to Lambert's works and their continuations. He also discusses a method published by Laplace in 1780.

{\bf 1798.} \citeauthor{laplace} publishes his {\it M{\'e}canique c{\'e}leste}. In \S 27 of the second book, he gives a proof of Lambert's theorem which is similar to Lagrange's proof ${\cal P}_4$. He concludes with three formulas, the first for the elliptic case, calling attention to the choices of arcs, the second for the parabolic case, being formula (\ref{Euler}) where the choices of signs are characterized, the third for the hyperbolic case. His discussion of signs includes a discussion of the extended rectilinear solutions. He republishes his orbit determination method, which does not use Lambert's theorem.

{\bf 1809.} \citeauthor{gauss} publishes his {\it Theoria motus}, a book on orbit determination. In \S 106 he gives the correct attribution of (\ref{Euler}): ``This formula appears to have been first discovered, for the parabola, by the illustrious Euler, (Miscell.\ Berolin, T.\ VII.\ p.\ 20,) who 
nevertheless subsequently neglected it, and did not extend it to the ellipse and 
hyperbola: they are mistaken, therefore, who attribute the formula to the illustrious Lambert, although the merit cannot be denied this geometer, of having 
independently obtained this expression when buried in oblivion, and of having 
extended it to the remaining conic sections. Although this subject is treated by 
several geometers, still the careful reader will acknowledge that the following 
explanation is not superfluous. We begin with the elliptic motion.'' Gauss gives a proof of Lambert's theorem of style ${\cal P}_4$. He insists on a remaining ambiguity of sign, which he explains by the existence of two ellipses: the second focus is constructed as the intersection of two circles, giving two possible positions. In \S 108 he discusses the limiting process to get formula (\ref{Euler}) from the elliptic case, but decides to give a proof of style ${\cal P}_1$, discussing the signs. In \S 109 he gives a proof and formulas for the hyperbolic case. He advertises the same expansion as Lambert, which is valid for the three conic sections, as being suitable if the orbit is nearly parabolic.
This short account of Lambert's theorem ends Gauss's study (\S 84--109) of what is called today the Lambert problem (see our remark~7). Gauss does not use Lambert's theorem to solve this problem, since he prefers the two methods he presents in \S 85--87 and \S 88--105.

{\bf 1815.} \citeauthor{lagrange7}, in the second, posthumous, edition of his {\it M{\'e}canique analytique}, section VII, \S 25, gives a proof of style ${\cal P}_1$ of formula (\ref{Euler}), cites Euler's article in \S 26, and presents a method for orbit determination which uses (\ref{Euler}). In \S 84 he presents briefly his proof ${\cal P}_3$. He shows how his final integral formula applies to the three kinds of conic sections, and advertises the same expansion as Lambert and Gauss. 

{\bf 1820.} \citeauthor{legendre} gives a proof of (\ref{Euler}) on page 7 of a book where he also recalls that he published the least-squares method four years before Gauss.

{\bf 1831.} \citeauthor{encke} gives a proof of (\ref{Euler}) in an article explaining Olbers's method. He introduces it as follows: ``Although this method was already carried to such a degree of perfection in the first memoir, that even the master-hand of the author of the {\it Theoria motus, \&c.}, made no essential alteration in it, but only some abbreviations, [...] Lambert's theorem is a main part of Olbers's method. The manner of solving it given by Olbers admitting of some abbreviations, I shall begin with explaining this little improvement.''

{\bf 1834.} \citeauthor{hamilton1} studies the properties of what he calls the characteristic function, namely, the integral of $2T dt$, where $T$ is the kinetic energy. This is $w$ of our formula (\ref{w}) in the case of a point particle. In \S 15 he shows that $w$ on elliptic arcs depends on $\|\A\B\|$, $\|\O\A\|+\|\O\B\|$ and $H$. Together with the relation $\delta w/\delta H=t$, on which he insists in \S 2, this gives a new proof of Lambert's theorem. His method to deal with $w$ has common features with Lagrange's proof ${\cal P}_4$. See e.g.\ his equation (108).

{\bf 1837.} \citeauthor{jacobi1} (\S 7) presents Hamilton's formulas in another order. He uses Lambert's theorem to deduce the trigonometrical expression of $w$ that Hamilton used to deduce Lambert's theorem. He insists on the analogy of the expressions of $t$ and $w$ (see Tait, 1866). He deduces from the expression of $w$ elegant formulas for the initial and final velocity vectors as $v_\A=k+\rho \varepsilon_\A$ and $v_\B=k-\rho \varepsilon_\B$ respectively, where $k$ is a vector along the chord, $\varepsilon_\A$ and $\varepsilon_\B$ are unit radial vectors, and $\rho$ is a number.  See 1866, 1888, 1961 and Sect.\thinspace\ref{sect10.2}. He checks that $w$ satisfies the  Hamilton-Jacobi equation.

{\bf 1837.} \citeauthor{chasles} (IV, \S 37) opposes again, after Lagrange, analysis and geometry: ``Le c{\'e}l{\`e}bre Lambert, autre Leibnitz par l'universalit{\'e} et la profondeur de ses connaissances, doit {\^e}tre plac{\'e} au nombre des math{\'e}maticiens qui, dans un temps o{\`u} les prodiges de l'analyse occupaient tous les esprits, ont conserv{\'e} la connaissance et le go{\^u}t de la G{\'e}om{\'e}trie et ont su en faire les plus savantes applications. [...] Ces consid{\'e}rations g{\'e}om{\'e}triques sont simples, et cependant elles ont suffi pour conduire Lambert au th{\'e}or{\`e}me le plus important de la th{\'e}orie des com{\`e}tes, dont les d{\'e}monstrations qu'on en a donn{\'e}es depuis par la voie du calcul ont exig{\'e} toutes les ressources de l'analyse la plus relev{\'e}e\footnote{The celebrated Lambert, another Leibniz by the universality and depth of his knowledge, must be placed among the number of mathematicians who, in a time when the miracles of analysis occupied all minds, retained the knowledge and the taste of Geometry and understood how to make the most savant applications of it. [...] These geometric considerations are simple, and yet they sufficed to lead Lambert to the most important theorem of the theory of comets, whose proofs given later by others using calculations required all the resources of the most exalted analysis.}.''

{\bf 1847.} \citeauthor{hamilton3} states a ``Theorem of hodographic isochronism: If two circular hodo\-graphs, having a common chord, which passes through or tends towards a common centre of force, be cut perpendicularly by a third circle, the times of hodographically describing the intercepted arcs will be equal.'' We can rephrase the statement by using our Definition~4: If a circle ${\cal C}$ cuts orthogonally two hodographs ${\cal H}_1$ and ${\cal H}_2$ of two Keplerian orbits of the same Lambert cycle, its center is on $\O x$. The arcs cut  on ${\cal H}_1$ and ${\cal H}_2$ are described in the same time.

{\bf 1862.} \citeauthor{cayley1} gives a description of Lambert's original results and of the Lambert cycle, and a computational proof of Lambert's theorem similar to ${\cal P}_4$, which is guided by Lambert's constructions.

{\bf 1866.} In \citeauthor{jacobi2}'s famous book, which is a course he gave in K\"onigsberg in the winter 1842--43, edited from notes by Borchardt,  lecture 25 is devoted to Lambert's theorem and its proof. Jacobi separates the Hamilton-Jacobi equation in elliptic coordinates, with a focus at $\O$ and another at the initial point $\A$. He slightly changes the presentation of his formulas for the initial and final velocities (see 1837). He shows how the separation produces elliptic integrals if the foci are $\O$ and an arbitrary point, even if this arbitrary point is a second fixed center. Except for the introduction of Hamilton's characteristic function, the proof follows ${\cal P}_3$, by reversing the order of generality. Lagrange is not cited for his proofs of Lambert's theorem, but only for his related article \citeyearpar{lagrange1}. See our Sect.\thinspace\ref{sect10.1}.

{\bf 1866.} \citeauthor{sylvester1} sets the semimajor axis of ellipses equal to 1 and proves by direct computation of the Jacobian that $\|\O\A\|+\|\O\B\|$, $\|\A\B\|$ and $\Delta t$ are functionally dependent. This is a proof of Lambert's theorem which he considers to be close to Lagrange's proof ${\cal P}_4$. But he actually removes part of Lagrange's computation, replacing it by the simpler-minded computation of the Jacobian. Then he takes the eccentricity $e$ as a parameter of what we call a Lambert cycle. He states that $\Delta t$ does not depend on $e$ and evaluates $\Delta t$ at $e=1$.

{\bf 1866.} \citeauthor{sylvester2} presents his previous proof with these words: ``Notwithstanding this plethora of demonstrations I venture to add a seventh, the simplest, briefest, and most natural of all''. He reacts to Lagrange's and Chasles's arguments about the advantage of geometry: ``In the nature of things such advantage can never be otherwise than temporary. Geometry may sometimes appear to take the lead of analysis, but in fact precedes it only as a servant goes before his master to clear the path and light him on the way. The interval between the two is as wide as between empiricism and science, as between the understanding and the reason; or as between the finite and the infinite''. He proves the hyperbolic and parabolic cases of Lambert's theorem as he had proved the elliptic case, and continues as described in the long title of his paper.

{\bf 1866.} \citeauthor{hamilton4}, at article 419 of his posthumous book {\it Elements of Quaternions}, proves his theorem of hodographic isochronism (see 1847), and deduces Lambert's theorem from it. He then gives a proof, using variations, quaternions and hodographs, of a ``new form of Lambert's Theorem'': the principal function from $\A$ to $\B$, and consequently the energy $H$, depend on $\|\O\A\|+\|\O\B\|$, $\|\A\B\|$, the elapsed time $\Delta t$ and the mass $m$ of the attracting body; the characteristic function, and consequently  $\Delta t$, depend on $\|\O\A\|+\|\O\B\|$, $\|\A\B\|$, $H$ and $m$.

{\bf 1866.} \citeauthor{tait1} (\citeyear{tait1} or \citeyear{tait2}, p.\ 163), interprets the analogy between time and characteristic function in \citet{hamilton1}: ``while the time is proportional to the area described about one focus, the action is proportional to that described about the other.''

{\bf 1869.} \citeauthor{cayley2} resolves the ambiguity of sign pointed out by Gauss with a geometrical criterion. One should ask if the line passing through $\A$ and the second focus separates $\O$ from $\B$.

{\bf 1878.} \citeauthor{adams} publishes a proof of type ${\cal P}_4$ in the elliptic and hyperbolic cases, and then in the parabolic case by passing to the limit. He presents the same formulas as Gauss in a more transparent way. He notices  that three functions are expressed in terms of two quantities only, $u_\B-u_\A$ and $e\cos\bigl((u_\A+u_\B)/2\bigr)$. This recalls Sylvester's argument. This presentation is adopted in \citet{dziobek},  \citet{routh} and \citet{battin}.

{\bf 1884.} \citeauthor{ioukovsky} proposes a proof based on the variation of the characteristic function $w$. He uses the analogy between $t$ and $w$ pointed out in \citet{jacobi1} and \citet{tait1} instead of using Hamilton's relation $\delta w/\delta H=t$. He does not cite any authors except Euler and Lambert. As the proof involves the second focus, one should adapt it to each kind of conic section. 

{\bf 1884.} \citeauthor{catalan} presents a proof of Lambert's theorem of style ${\cal P}_4$, where he interprets each step with a geometrical construction. He gives some related geometrical statements, one of them being a construction,  from a general ellipse with a chord, of what we call the vertical branch $\Sigma$ (see 1761, 1862), others being new.

{\bf 1888.} \citeauthor{dziobek}'s book gives a short proof using Adams's argument and a proof inspired by Hamilton and Jacobi. He comments: ``For a long time, the proposition was regarded as a curiosity. Its true source was shown by the investigations of Hamilton and Jacobi.'' He advertises Jacobi's expression of $v_\A$ and $v_\B$ and writes: ``no one would have succeeded {\it a priori} in getting the notable equations [of $v_\A$ and $v_\B$] from those \S 1.'' We will comment on his words in Sect.\thinspace\ref{sect10.2}.

{\bf 1901.} \citeauthor{bourget} complains that \citet{jacobi2} does not cite Lagrange's proof ${\cal P}_6$. He generalizes the main identity in ${\cal P}_6$.

{\bf 1941.} \citeauthor{wintner}'s book, \S 247-248, introduces his precise presentation in this way: `A proof of Lambert's theorem can be obtained by an application of the theorem of Gauss-Bonnet on the surface of revolution ${\bf S}_h$ of \S 244. However, the proof is shorter if use is made of the ``Beltrami-Hilbert integral'' or the ``isoenergetic action $W$'' not via ${\bf S}_h$ but in a more direct manner, as follows.' His historical note on p.\ 422 compares the lengths of various proofs.

{\bf 1961.} \citeauthor{godal} presents as does Jacobi in 1837 the initial and final velocity vectors as $k+\rho \varepsilon_\A$ and $k-\rho \varepsilon_\B$ respectively. He notices that $\rho\|k\|$ depends only on $\A$ and $\B$, not on the orbit. See Sect.\thinspace\ref{sect10.2}.

{\bf 1966.} \citeauthor{levine}, while developing a method of orbital navigation,  considers the ray from the center $\O$ to a point where the velocity of the spacecraft is parallel to a chord $\A\B$ of the orbit. He notices that the angle from $\A\B$ to such a ray depends only on $\A$ and $\B$, not on the orbit.

{\bf 1976.} \citeauthor{correas} insists on proving Lambert's theorem for the three kinds of conic sections in a single argument, and does it by giving a single formula for $\Delta t$ by means of Stumpff's functions.

{\bf 1983.} \citeauthor{souriau} proposes (see p.\ 376) a new proof of Lambert's theorem, which uses a collection of remarkable and elegant formulas about the Kepler problem.

{\bf 2002.} \citeauthor{marchal} presents several formulas and proofs of style ${\cal P}_4$ which include, as \citet{jacobi1} does, formulas for the action. Remarkable inequalities are deduced and used to estimate the minimizers of the action in the $n$-body problem.

{\bf 2016.} Authors such as \citeauthor{LinetTeyssandier} show us that Lambert's theorem may still be rediscovered by skillful calculators. They consider the gravitational influence of a spherically symmetric body on the propagation of light within the weak-field, linear approximation of general relativity. Their formula (39) is typically ``Lambertian''.

\section{Final comments}
Many proofs were proposed after Lambert's proof in 1761. Such a ``plethora of demonstrations'', in Sylvester's words, gives the impression of a chronic dissatisfaction. After Lambert's publication, which was found to be obscure, most attempts were ``analytical''. The  geometrical arguments of Lambert remained essentially untouched, being only described in few words by \citet{lagrange4} and in a short note by \citet{cayley1}. The fact that two unparametrized arcs belonging to the same Lambert cycle correspond to each other through an affine transformation of the plane has apparently never been stated.

\subsection{Comments on the first proof by Lagrange}
\label{sect10.1}
The elliptic coordinates $\sigma$ and $\tau$ of a point moving in a plane attracted by two Newtonian fixed centers are two elliptic functions of a common parameter, which is not the time, while the time parameter is expressed as an elliptic integral in $\sigma$ minus an elliptic integral in $\tau$. By contrast, when there is only one fixed center, the motion is Keplerian. The analytic expression of the position of the moving point is simpler and does not involve any elliptic function. Consequently, when one of the two fixed centers has zero mass, simplifications should occur when combining the elliptic functions. Observing these simplifications, \citeauthor{lagrange2} obtained a proof of Lambert's theorem which was published in 1773.  Indeed, when one of the masses is zero, $\sigma$ and $\tau$ are expressed by the {\it same} elliptic function, with a constant shift of the common parameter, and the time is expressed as an elliptic integral in $\sigma$ minus the {\it same} elliptic integral in $\tau$. A strange simplification occurs when subtracting. The formulas appear as (M), (N) and (T) in \citet{lagrange2}, and in \S 10 of \citet{lagrange4}. The deduction of (T) uses a method explained in \citet{lagrange1} without  reference to the mechanical problem, but with a reference to previous works by Euler about elliptic integrals. \citet{euler3} commented on Lagrange's recent works in his last letter to him. These works became classical in the theory of elliptic functions, about the addition theorem (see the notes by the editors of \citealt{euler3}, and \citealt{houzel}, p.\ 89).  They concern the Keplerian motion expressed in elliptic coordinates rather than Lambert's theorem. As we said, \citet{jacobi2} uses \citet{lagrange1}. \citet{sylvester2} advertises Lagrange's identities without mentioning Euler or Lexell. \citet{bourget} cites works by Euler, Raffy, Fagnano, Graves and Chasles.

\subsection{Comments on our minimal proof}
\label{sect10.2}
\citet{dziobek} claims that Hamilton and Jacobi found the ``true source'' of Lambert's theorem. He is not convincing: of the two proofs in his book, the short one does not involve such a ``source'', while the long one does. If Hamilton himself  was convinced he got the ``true source'' in 1834, he would not have published other proofs based on different ideas. What Hamilton indeed shared with his contemporaries is an obsession with Lambert's theorem. Uncovering deep features of dynamics and geometry, namely, the properties of the characteristic function, the circular hodograph of the Keplerian motion and the quaternion algebra, he successively used them to produce new demonstrations. Jacobi  does not appear to be convinced in 1837 that Hamilton got the ``true source'', and the key to Jacobi's second proof is the elliptic system of coordinates rather than the characteristic function.

The simplicity of our minimal proof supports Dziobek's opinion about the ``true source'' and at the same time contradicts his words ``no one would have succeeded'' (see 1888). Lemma 4 is remarkable. The direction of $v_\B-v_\A$ does not depend on the choice of the conic section passing through $\A$ and $\B$.  If $\varepsilon_\A+\varepsilon_\B\neq 0$, there is a $\rho$ such that $v_\A-v_\B=\rho (\varepsilon_\A+\varepsilon_\B)$. We set $k=v_\A-\rho \varepsilon_\A=v_\B+\rho \varepsilon_\B$ and get Jacobi's expressions $v_\A=k+\rho \varepsilon_\A$, $v_\B=k-\rho \varepsilon_\B$, where it just remains to express $k$ and $\rho$, if needed. Interestingly, Dziobek refers to his \S 1 as not giving this key lemma, but this first section of his excellent book does present the eccentricity vector in (17a), in a new and deep way, and does use it to compute velocities, in his proof of the circularity of the hodograph.

\medskip\noindent{\bf Remark 15.} \citet{hamilton1} obtained from the expression of the characteristic function ``the following curious, but not novel property, of the ellipse'', which is republished in \citeyear{hamilton4}, just after the ``new form of Lambert's theorem'', as ``this known theorem: that {\it if two tangents $({\rm QP},{\rm QP}')$ to a conic section be drawn from any common point $(\Q)$, they subtend equal angles at a focus $(\O)$}, whatever the special form of the conic may be''. One should understand that the equal angles are ${\rm QOP}$ and ${\rm QOP}'$. 

\vspace{0.5cm}
\centerline{\includegraphics[width=40mm]{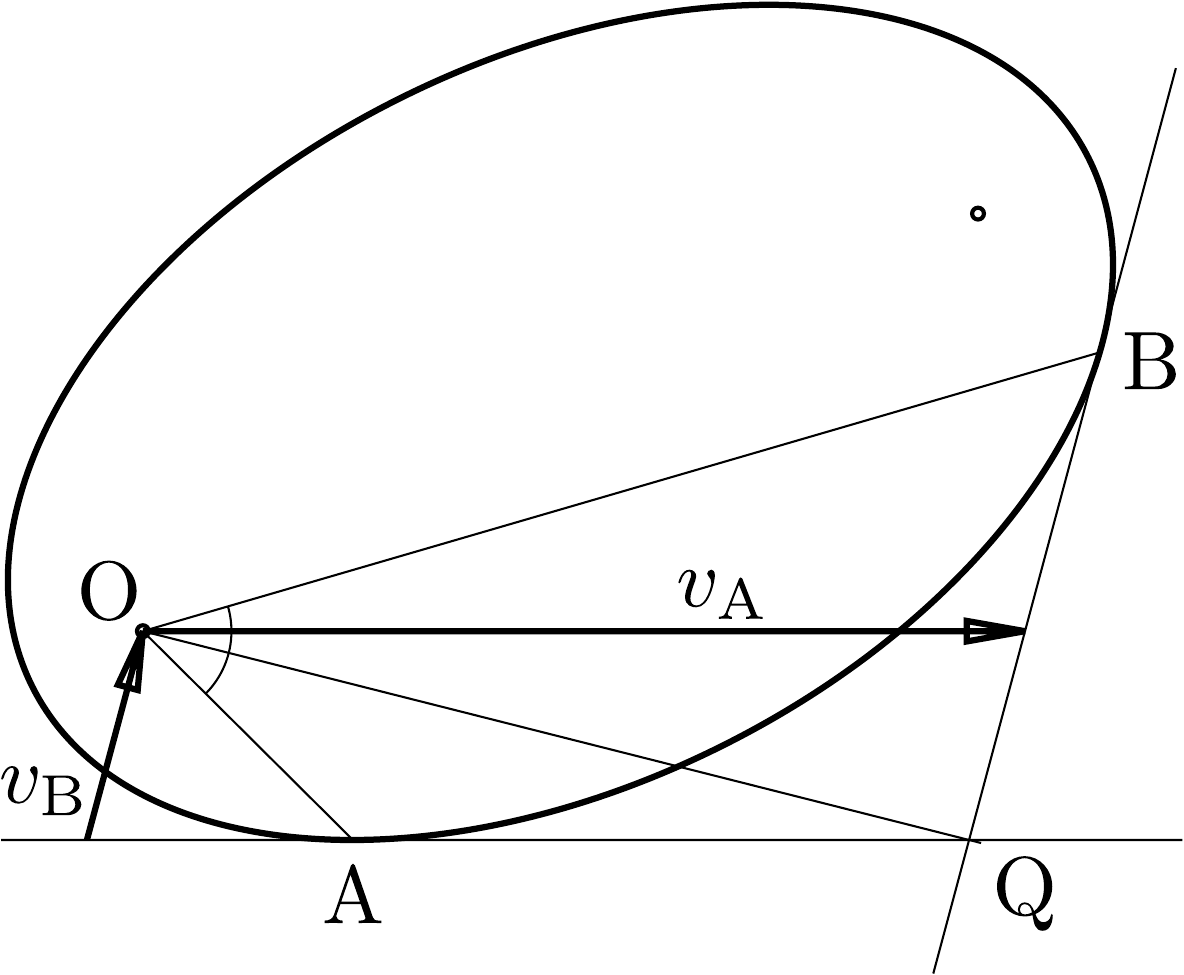}}
\centerline{\sl Fig.\thinspace11. Lemma 12, or Hamilton's statement with $\A={\rm P}$, $\B={\rm P}'$.}
\vspace{0.5cm}

The same property  appears in \citet{ioukovsky}, now as an argument used to prove Lambert's theorem. According to \citet{berger}, 17.2.1.6, this property is one of Poncelet's ``small theorems''. \citet{poncelet}, p.\  265, states this property and the fact that the external bisector of ${\rm POP}'$ meets the chord on the directrix, but he gives credit to De Lahire and l'H{\^o}pital. The earliest statement we know of the ``curious property'' belongs to \citet{hire}, book 8, Proposition 24, p.\ 190.

Consider Fig.\thinspace11. The pair of velocity vectors at $\A$ and $\B$ should be proportional to the represented pair since $\O\A\wedge v_\A=\O\B\wedge v_\B$. But $v_\A-v_\B=\O\Q$. De La Hire's property is thus reduced to Lemma 4, which can be restated as:

\medskip\noindent{\bf Lemma 12.} Consider two positions $\A$ and $\B$ on a Keplerian orbit in a plane with origin the fixed center $\O$. Let  $v_\A$ and $v_\B$ be the velocity vectors at these positions. The interior bisector line of the angle $\A\O\B$ is directed along $v_\B-v_\A$ and passes through the intersection $\Q$ of the respective tangents at $\A$ and $\B$.

\subsection{Comments on our second proof and the question in the title}
This constructive proof  improves Lambert's original proof ${\cal P}_2$ and Lagrange's proof ${\cal P}_4$, by pointing out the affine transformations and Fig.\thinspace8, and by getting the three kinds of conic sections in a single computation. Note that ${\cal P}_4$ has longer computations than our proof only for the elliptic case.

Our fundamental identity (\ref{o11}) uses in a non-intuitive way the most typical operation of Algebra, the ``al-jabr'' operation, which consists in translating a term from the left-hand side to the right-hand side of an equation. We were not able to find a purely geometrical argument of comparable simplicity. In all other attempts, the parabolic case required a special treatment. We presented our Proposition 7 as a partial success in an attempt to include the parabolas in a geometrical statement related to Lambert's theorem.

The solution of the Kepler problem is pure geometry as far as the time is ignored. The time parametrizes the solutions transcendently. Lambert's theorem gives a geometric property of the time. We cited Lagrange, Chasles and Sylvester discussing the question: should this Theorem be proved by geometry or analysis?

We propose a related question: does Theorem 1 belong to geometry or to dynamics?  In all the attempts to remove the time parameter, exceptions concerning the parabolic and the rectilinear orbits appear, which complicates the statement. Lambert's theorem generates theorems on conic sections, but conversely we are not able to present it as a simple corollary of a theorem on conic sections. This suggests that a chapter of elementary geometry finds its source in classical dynamics.

\medskip\noindent{\it Acknowledgements.} Thanks to Alain Chenciner for many helpful suggestions, to Richard Montgomery for his encouraging help, to H.\ Scott Dumas for the translations in the footnotes and his comments, to R{\'e}mi Bourgeois, Pierre Teyssandier and Christian Velpry for criticizing my drafts, to Zhao Lei for many remarks. Thanks also to the students of AIMS S{\'e}n{\'e}gal for criticizing my oral presentation of the constructive proof, to Niccol{\`o} Guicciardini and Jes{\'u}s Palac{\'\i}an for valuable references, and to the editors and reviewers for their suggestions.

\end{document}